\documentclass[aop,preprint]{imsart}

\RequirePackage{amsthm,amsmath,amsfonts,amssymb}
\RequirePackage[numbers]{natbib}
\RequirePackage[colorlinks,citecolor=blue,urlcolor=blue]{hyperref}

\usepackage{booktabs}

\usepackage{enumitem}
\usepackage{amsmath}
\usepackage{amsfonts}
\usepackage{amssymb}             
\usepackage{dsfont}

\usepackage{comment}
\usepackage{breakurl}

\usepackage{mathrsfs}
\usepackage{booktabs}
\usepackage{mathtools}

\usepackage{tikz}
\usepackage{amsthm}                
\usepackage{mhequ}
\usepackage{hyperref}

\startlocaldefs

\hypersetup{pdfstartview=}

\addtocontents{toc}{\protect\hypersetup{hidelinks}}

\DeclareSymbolFont{sfoperators}{OT1}{ptm}{m}{n}
\DeclareSymbolFontAlphabet{\mathsf}{sfoperators}

\makeatletter
\def\operator@font{\mathgroup\symsfoperators}
\makeatother
\newcommand{\eqdef}{\stackrel{\mbox{\tiny def}}{=}}

\numberwithin{equation}{section}

\newcommand{\R}{\mathbb{R}}

\newcommand{\N}{\mathbb{N}}

\newcommand{\cA}{\mathcal{A}}

\newcommand{\CC}{\mathcal{C}}

\newcommand{\cS}{\mathcal{S}}

\newcommand{\dd}{\mathrm{d}}      

\newcommand{\eps}{\varepsilon}

\newcommand{\Di}{\mathrm{Diag}}

\newcommand{\oDi}{\mathrm{off}}

\colorlet{darkblue}{blue!90!black}
\colorlet{darkred}{red!90!black}
\colorlet{darkgreen}{green!50!black}
\colorlet{darkyellow}{yellow!90!black}

\newcommand{\gena}{\cA}
\newcommand{\genap}{\cA_+}

\newcommand{\Op}{\mathcal{H}}

\newcommand{\Ll}{\mathrm{L}}
\newcommand{\LB}{\mathrm{LB}}
\newcommand{\UB}{\mathrm{UB}}

\newtheorem{thm}{Theorem}[section]

\newtheorem{lem}[thm]{Lemma}

\newtheorem{assumption}[thm]{Assumption}
\theoremstyle{remark}
\newtheorem{remark}[thm]{Remark}

\endlocaldefs

\begin{document}

\begin{frontmatter}
\title{$\lowercase{\sqrt{\log t}}$ - superdiffusivity for a Brownian particle in the curl of the 2d GFF}
\runtitle{Superdiffusivity for a Brownian particle in the curl of the GFF}

\begin{aug}
\author[A]{\fnms{Giuseppe} \snm{Cannizzaro}\ead[label=e1]{giuseppe.cannizzaro@warwick.ac.uk}},
\author[B]{\fnms{Levi} \snm{Haunschmid-Sibitz}\ead[label=e2,mark]{levi.haunschmid@tuwien.ac.at}}
\and
\author[B]{\fnms{Fabio} \snm{Toninelli}\ead[label=e3,mark]{fabio.toninelli@tuwien.ac.at}}
\address[A]{University of Warwick, UK, \printead{e1}}

\address[B]{Technical University of Vienna, Austria, \printead{e2,e3}}
\end{aug}

\begin{abstract}
The present work is devoted to the study of the large time behaviour of a  critical Brownian
diffusion in two dimensions, whose drift is divergence-free, ergodic and
given by the curl of the 2-dimensional Gaussian Free Field. We prove the conjecture,
made in [B. T\'oth, B. Valk\'o, J. Stat. Phys., 2012], according to which the diffusion
coefficient $D(t)$ diverges as $\sqrt{\log t}$ for $t\to\infty$.
Starting from the fundamental work
  by Alder and Wainwright [B. Alder, T. Wainwright, Phys. Rev. Lett. 1967],
  logarithmically superdiffusive behaviour has
  been predicted to occur for a wide variety of out-of-equilibrium
  systems in the critical spatial dimension $d=2$. Examples include the diffusion of a tracer particle in a fluid,
  self-repelling polymers and random walks, Brownian particles
  in divergence-free random environments, and, more recently,
  the 2-dimensional critical Anisotropic KPZ equation.
  Even if in all of these cases it is expected that $D(t)\sim\sqrt{\log t}$,
  to the best of the authors' knowledge, this is the first instance in which
  such precise asymptotics is rigorously established.
\end{abstract}



\end{frontmatter}

\section{Introduction}

\label{sec:intro}
In the present work, we study the motion of a Brownian particle in $\mathbb R^2$,
subject to a random, time-independent
drift  $\omega$ given by the curl of the two-dimensional Gaussian Free Field ($2d$ GFF).
Namely, we look at the SDE which is (formally) given by
\begin{equ}
  \label{eq:SDE}
  \dd X(t)=\omega(X(t))\dd t +  \dd B(t), \qquad X(0)=0
\end{equ}
where $B(t)$ is a standard two-dimensional Brownian motion and
\[
x\mapsto \omega(x)=(\omega_1(x),\omega_2(x))
\]
is defined as
\begin{equ}
  \label{eq:omegaxi}
x=(x_1,x_2)\mapsto \omega(x)=(\partial_{x_2}\xi(x),-\partial_{x_1}\xi(x))\,,
\end{equ}
with $\xi$ the $2d$ GFF. As written,~\eqref{eq:SDE} is ill-posed due to the singularity
of the drift $\omega$. In fact, not only classical stochastic analytical tools would
fail in characterising (even) its law but it would also be {\it critical} for the recent techniques
established in~\cite{CC,DD} as its spatial regularity is way below the threshold identified
therein\footnote{Formally, the $2d$ GFF is in $\CC^{\alpha}$, $\alpha<0$, the latter being the
space of H\"older distributions with regularity $\alpha$ (see~\cite{CC} for the definition), so that
$\omega\in\CC^{\alpha-1}$. In the aforementioned works, the threshold regularity is $-2/3$ so that~\eqref{eq:SDE} falls indeed out of their scope.}.
Nevertheless, we are interested in its large time behaviour and hence
we regularise $\xi$ by convolving it with a $C^\infty$ bump function (see
Section \ref{sec:themo} for details), so that $\omega$ is well-defined
pointwise and smooth. Note that the vector field $\omega$ is
everywhere orthogonal to the gradient of the field $\xi$, and
therefore parallel to its level lines. As a consequence, the particle is
subject to two very different mechanisms: the drift tends to push the
motion \emph{along the level lines of the GFF}, while the Brownian
noise tends to make it diffuse \emph{isotropically}. Our main theorem is a
{\it sharp} superdiffusivity result: the mean square displacement
$\mathbf E[|X(t)|^2]$ (under the joint law of the Brownian noise and of
the random drift) is of order $t \sqrt{\log t}$ for $t\to\infty$, up
to multiplicative loglog corrections. This proves a conjecture of B. T\'oth and
B. Valk\'o \cite{Toth2012Superdiffusive2} and, in a broader
perspective, it is the first proof of the $\sqrt{\log t}$-{\it superdiffusivity phenomenon}
conjectured to occur in a large class of
(self-)interacting diffusive systems in dimension $d=2$ (see the
discussion below).

To put the model and the result into context, let us observe first
that the vector field $\omega$ is divergence-free and that its law
is translation-invariant and ergodic. Brownian diffusions in ergodic,
divergence-free vector fields have been introduced in the physics and
mathematics literature as a (toy) model for a tracer particle evolving
in an incompressible turbulent flow.  If the energy spectrum of the
vector field (i.e. the Fourier transform $e(p)$ of the trace of the
covariance matrix
$R(x-y)=\{\mathbb E(\omega_a(x)\omega_b(y))\}_{a,b\le d}$, with $d$
the space dimension) satisfies the integrability condition\footnote{The integral in the l.h.s. is known as ``P\'eclet number''
\cite{komorowski2012fluctuations}.}
\begin{equ}
  \label{eq:integra}
  \int_{\mathbb R^d} \frac{e(p)}{|p|^2}\dd p<\infty,
\end{equ}
 the behaviour of the particle is
known to be diffusive on large scales
\cite{komorowski2001homogenization,komorowski2012fluctuations} (see also \cite{kozma2017central,toth2018quenched}
for analogous results obtained, via different methods,
in the discrete setting of random walks in divergence-free random environments). In the
robustly superdiffusive case, where the integral in
\eqref{eq:integra} has a power-law divergence for small $p$, it turns out that
$\mathbf E[|X(t)|^2]$ grows like $t^{\nu}$ for some $\nu>1$   \cite{komorowski2002superdiffusive}. The case
under consideration in this work instead, where $d=2$ and $\omega$ is
the curl of the GFF, is precisely at the boundary between the
diffusive and the super-diffusive case: $e(p)$ is essentially constant  for $p$ small, the integral
\eqref{eq:integra} diverges logarithmically at small momenta and
logarithmic corrections to diffusivity are expected.

Logarithmic corrections to diffusivity in two-dimensional
out-of-equilibrium systems have a long history.
The seminal works~\cite{alder1967velocity,wainwright1971decay} of Alder and Wainwright lead the way,
in that they predicted
that the velocity auto-correlation of a tracer particle diffusing in a
fluid behaves like $t^{-d/2}$ in dimension $d\ge 3$ and like
$1/(t\sqrt{\log t})$ in the critical dimension $d=2$. This translates
into the fact that, in two dimensions, the mean square displacement
$\mathbf E[|X(t)|^2]$ of the particle should grow like $t D(t)$ with
\begin{equ}[eq:conjecture]
D(t)\approx \sqrt{\log t} \qquad \text{as} \quad t\to\infty.
\end{equ}
The quantity $D(t)$ takes the name of (bulk) diffusion coefficient.
The same prediction was obtained by Forster, Nelson and Stephen
\cite{forster1977large} via Renormalization Group methods. Subsequently, anomalous
logarithmic corrections as in \eqref{eq:conjecture} were conjectured to occur for several other
two-dimensional (self-)interacting diffusions, including
self-repelling random walks and Brownian polymers\footnote{There has been some controversy in the physics literature as to the value of the exponent $\zeta$ of the logarithm in \eqref{eq:conjecture} for self-repelling random walks. The values $\zeta=0.4$ and $ \zeta=1$  have been proposed \cite{amit1983asymptotic,obukhov1983renormalisation,peliti1987random}, in addition to the $\zeta=1/2$ prediction \cite{Toth2012Superdiffusive2} based on the Alder-Wainwright argument.} \cite{amit1983asymptotic,obukhov1983renormalisation,peliti1987random,Toth2012Superdiffusive2}, lattice gas models \cite{landim2005superdiffusivity},
the diffusion \eqref{eq:SDE} in the curl of the $2d$ GFF \cite{Toth2012Superdiffusive2} and, more
recently, the two-dimensional Anisotropic KPZ equation ($2d$ AKPZ)
\cite{Cannizzaro2020TheSuperdiffusivity}. We emphasize that in all
of these cases, it is known  or conjectured that the analogous models  behave
diffusively ($D(t)\sim 1$) in dimension $d\ge 3$
(see for instance \cite{Horvath2012Diffusive3} for the self-interacting random walks and Brownian polymers).

From a rigorous viewpoint, results available so far fall short
of the conjecture \eqref{eq:conjecture}. Until recently, the best estimates obtained can be summarised
into bounds of the form
 \begin{equ}[eq:previous]
   \log\log t\lesssim D(t)\lesssim\log t
 \end{equ}
 (see \cite{Toth2012Superdiffusive2} for  $2d$ self-repelling Brownian polymers and for the SDE \eqref{eq:SDE}, and \cite{landim2005superdiffusivity} for two-dimensional lattice fluids). More recently, two of the authors together with D. Erhard proved in~\cite{Cannizzaro2020TheSuperdiffusivity} that, for the $2d$ AKPZ equation, one has
 \begin{equ}
   \label{eq:previous2}
 (\log t)^a\lesssim D(t)\lesssim (\log t)^{1-a}
 \end{equ}
 for some sufficiently small $a>0$; after the present work was
 completed, in a second version of
 \cite{Cannizzaro2020TheSuperdiffusivity} the result for the $2d$ AKPZ equation has been also improved
 to $a=1/2$.  (All the above cited results have been shown in the
 sense of Laplace transform.)

For the SDE \eqref{eq:SDE}  under consideration in the present work,
we establish for the first time the conjectured behaviour
\eqref{eq:conjecture}, up to corrections that are polynomial in
$\log \log t$ (see Theorem~\ref{th:main} below). The result holds again in the sense of Laplace
transform - see, however, Remark \ref{rem:realtime} for its
implications in real time.
\newline

Our argument is based on an iterative analysis of the resolvent of the
generator of the Markov process given by the environment seen from the
particle (see~\eqref{eq:seen} below). This is inspired by the method employed by
H.-T. Yau \cite{Yau2004logtProcess} to prove $(\log t)^{2/3}$ corrections to the
diffusivity of the two-dimensional Asymmetric Simple Exclusion
Process ($2d$ ASEP) and, more closely, by the techniques developed in
\cite{Cannizzaro2020TheSuperdiffusivity} to determine~\eqref{eq:previous2}
for the $2d$ AKPZ equation. Note that the exponent $2/3$ of the
logarithmic corrections of $2d$ ASEP is different from the exponent
$1/2$ in \eqref{eq:conjecture}, reflecting the fact that the two
models belong to two different universality classes, as emphasized
already in \cite{landim2005superdiffusivity,Toth2012Superdiffusive2}. From a technical point of view, a crucial
difference between the two models is that for $2d$ ASEP the iterative
method in \cite{Yau2004logtProcess} provides, at each step $k$ of the recursion,
upper/lower bounds for $D(t)$ of the form $(\log t)^{\nu_k}$, with
$\nu_k$ converging exponentially fast to $2/3$ as $k\to\infty$. In our
case, on the other hand, at step $k$ the method naturally provides
lower (resp. upper) bounds of order $(\log \log t)^k/k!$ (resp.
$k! \log t/(\log \log t)^k$) and we have to run the iteration for a
number of steps of order $k=k(t)\approx \log\log t$ (instead of
$k(t)\approx\log\log \log t$ as in \cite{Yau2004logtProcess}) to reach the final
result. As a consequence, in contrast with \cite{Yau2004logtProcess}, we cannot
afford to lose a multiplicative constant at each step of the iteration
(such multiplicative constants are responsible for the sub-optimal
result \eqref{eq:previous2} in the first version of
\cite{Cannizzaro2020TheSuperdiffusivity}), and a much finer analysis
of the resolvent is needed.  Further, we get a significantly
sharper control of sub-leading corrections to $D(t)$ with respect to
$2d$ ASEP, namely, a multiplicative correction that is polynomial in
$\log \log t$ (see Remark \ref{rem:realtime}), to be compared with the
corrections of order $\exp( (\log \log\log t)^2)$ for $2d$ ASEP
\cite{Yau2004logtProcess}.

\subsection*{Organization of the article}
The rest of this work is organized as follows. In Section
\ref{sec:themo}, we rigorously define the model~\eqref{eq:SDE}
and state the main result. In Section~\ref{sec:prelim}, we introduce the main tools:
we recall the generator of the environment seen from the particle process and we describe the
space on which it acts. Section~\ref{sec:gen} is devoted to the analysis
of the generator and the derivation of the crucial recursive bounds, while in Section~\ref{sec:Proof},
the proof of the main result is given. At last, in Appendix~\ref{app}, we collect some technical estimates needed
in Sections~\ref{sec:gen} and~\ref{sec:Proof}.

\section{The model and main result}
\label{sec:themo}
The Brownian diffusion in the curl of the 2-dimensional Gaussian Free Field  is the stochastic process $t\mapsto X(t)\in\mathbb R^2$ given by the solution of the SDE \eqref{eq:SDE}
where $B(t)$ is a standard two-dimensional Brownian motion and
\[x\mapsto \omega(x)=(\omega_1(x),\omega_2(x))\] is a smooth,
divergence-free, random vector field on $\mathbb R^2$, given by the
curl of (a smoothened version of) the two-dimensional Gaussian free
field. To be more precise, let us introduce the following assumption which will be
in place throughout the paper.

\begin{assumption}\label{a:V}
Let $V\colon \R^2\to\R$ be a radially symmetric bump function, i.e. a function such that
there exists $U\colon \R^2\to\R$ which is a smooth function
in $C^\infty(\R^2)$, radially symmetric, decaying sufficiently (say, exponentially) fast at infinity and such that
$\int_{\mathbb R^2} U(x)\dd x=1$, for which
\begin{equ}[e:V]
V=U\ast U\,.
\end{equ}
\end{assumption}

Let $V$ satisfy Assumption~\ref{a:V} and $U$ be such that~\eqref{e:V} holds.
Then, $\omega$ is a centred Gaussian field whose law
$\mathbb P$ (and corresponding expectation $\mathbb E$), is defined as follows.
Let $\xi\colon\mathbb R^2\to \mathbb R$ be the two-dimensional Gaussian
Free Field convoluted with $U$, i.e. the centred Gaussian field with covariance
\begin{equ}
\mathbb E(\xi(x)\xi(y))=V * g(x-y)\,,\qquad \text{for all $x,\,y\in\R^2$}
\end{equ}
where $g(x)=-\log |x|$.
Then, we define $\omega $ as the curl of the scalar field $ \xi$, i.e. as in \eqref{eq:omegaxi},
which clearly satisfies for all $x,\,y\in\R^2$
\begin{equ}
  \label{eq:covxi}
  \mathbb E(\omega_\ell(x))=0, \quad \mathbb E(\omega_k(x)\omega_\ell(y))=-\tilde\partial_{k}\tilde\partial_{\ell}V*g(x-y), \quad k,\ell=1,2,
\end{equ}
for $\tilde\partial_{1}:=\partial_{x_2}$ and $\tilde\partial_2:=-\partial_{x_1}$.
Note that, while convolving the full-plane Gaussian free field with $U$ is a somewhat formal operation (because the field is only defined up to a constant), the derivatives of the convolved field and therefore $\omega$ are (pointwise) defined without any ambiguity and are smooth with respect to $x$.

It is well known \cite[Chapter 11]{komorowski2012fluctuations} that, since $\omega$
sampled from $\mathbb P$ is divergence-free, translation invariant and
ergodic, the law $\mathbb P$ is stationary for the Markov process of the environment seen
from the particle, i.e. the time-evolving field
$t\mapsto \omega(t,\cdot)$ given by
\begin{equ}
  \label{eq:seen}
  \omega(t,x)\eqdef\omega(X(t)+x)\,,\qquad x\in\R^2\,.
\end{equ}

\subsection{Main result}

Our main result is a sharp estimate on the super-diffusivity of the process $X$. For $t>0$, let
$  \mathbf E(|X(t)|^2)$ denote the mean square displacement of $X$ at time $t$ - the expectation being taken
with respect to the joint randomness of the vector field $\omega$ and of the Brownian noise $B$ in \eqref{eq:SDE}.
Let us remark that  $\mathbf E(X(t))=0$ because the law of the environment is symmetric  and $X(0)=0$.

Throughout the present article we will be working with its Laplace transform, given by
\begin{equ}
\label{eq:LaplaceT}
  D(\lambda)\eqdef\int_0^\infty e^{-\lambda t} \mathbf E(|X(t)|^2)\dd t, \quad \lambda>0.
\end{equ}
Note that if in~\eqref{eq:SDE} there were no drift, one would trivially have $X(t)=B(t)$
so that $ \mathbf E(|X(t)|^2)=2 t$ and $D(\lambda)=2 \lambda^{-2}$.
In \cite{Toth2012Superdiffusive2}, it was conjectured that
\begin{equ}
  \label{eq:conj}
   D(\lambda)\stackrel{\lambda\to0}\approx {\lambda^{-2}}\sqrt{|\log \lambda|},
\end{equ}
corresponding in real time to
\begin{equ}
  \mathbf E(|X(t)|^2)\stackrel{t\to\infty}\approx t\sqrt{\log t}
\end{equ}
(the diffusion coefficient mentioned in the abstract is $D(t)= t^{-1}\mathbf E(|X(t)|^2)$.
From a rigorous point of view, in the aforementioned work it was proved that
\begin{equ}
  C_1 \lambda^{-2}\log|\log\lambda|\le D(\lambda)\le C_2 \lambda^{-2}|\log \lambda|
\end{equ}
for some positive constants $C_1,C_2$, for sufficiently small $\lambda$.
In this work, we establish the conjecture \eqref{eq:conj} in full.
\begin{thm}
  \label{th:main}
  For every $\eps>0$ there exists constants $C_\pm(\eps)$ such that, for every $0<\lambda<1$,
  \begin{equ}
    \label{eq:main}
    C_-(\eps)\big(\log|\log\lambda|\big)^{-1-\eps}
    \le \lambda^2\frac{D(\lambda)}{\sqrt{|\log\lambda|}} \le  C_+(\eps)\big(\log|\log\lambda|\big)^{1+\eps}.
  \end{equ}
  where $D$ is defined according to~\eqref{eq:LaplaceT}.
\end{thm}
  The exponent $1+\eps$ in the sub-dominant corrections can presumably be improved by some additional technical work,
  but we do not pursue this here.
The constants $C_\pm(\eps)$ implicitly depend also on the choice of bump function $V$.
  \begin{remark}
    \label{rem:realtime}
    By a well-established argument (see~\cite{MR2477398}) the upper bound in~\eqref{eq:main}
    implies an upper bound for the diffusivity in real time of the form
    \begin{equ}
       \mathbf E(|X(t)|^2)\le  O\left(t \sqrt{\log t}\;(\log\log t)^{1+\eps}\right).
    \end{equ}
    Deducing a pointwise (in time) lower bound on
    $\mathbf E(|X(t)|^2)$ from the behaviour for $\lambda\to0$ of the
    Laplace transform is much more delicate. That said, one can easily
    get (applying for instance \cite[Theorem
    1.7.1]{bingham1989regular}) the following
    \begin{equ}
      \limsup_{t\to\infty}\frac{ \mathbf E(|X(t)|^2)}{t \sqrt{\log t}\;(\log\log t)^{-1-\eps}}>0\,.
    \end{equ}
\end{remark}

\section{Preliminaries}
\label{sec:prelim}
By rotation invariance, one has $\mathbf E(|X(t)|^2)= \mathbf E(X_1(t)^2+X_2(t)^2)= 2\mathbf E(X_1(t)^2)$,
and we write
\begin{equ}
  \label{eq:soluzione}
  X_1(t)=B_1(t)+\int_0^t \phi(\omega_s) \dd s=:B_1(t)+F_1(t)
\end{equ}
where $t\mapsto \omega_t$ is the environment seen from the particle (recall \eqref{eq:seen}), and
\begin{equ}
  \label{eq:phi}
\phi(\omega)\eqdef\omega_1(0)\,.
\end{equ}
Recall that $\mathbf E X_1(t)=0$. The first term on the r.h.s. of
\eqref{eq:soluzione} has variance $t$, so to prove~\eqref{eq:main}, it is sufficient to show
\begin{equ}
  \label{eq:main2}
  \frac{ C_-(\eps)}{\big(\log|\log\lambda|\big)^{1+\eps}}
     \le \frac{\lambda^2}{\sqrt{|\log\lambda|}}\int_0^\infty e^{-\lambda t} \mathbf E[F_1(t)^2]\dd t \le  C_+(\eps)\big(\log|\log\lambda|\big)^{1+\eps}\,.
\end{equ}

The starting point in the study of $\mathbf E[F_1(t)^2]$ is the
understanding of the environment process $t\mapsto \omega_t$. As argued in
\cite{Toth2012Superdiffusive2}, this is a Markov process, whose generator will be denoted by
$\mathcal G$, on the Fr\'echet space of $C^\infty$, divergence-free
two-dimensional vector fields with derivatives growing slower than any
power at infinity.
As the field is stationary, ergodic and divergence-free,
the probability measure $\mathbb P$ is
stationary for the environment process  \cite[Chapter 11]{komorowski2012fluctuations}.
This ensures that, as in \cite[Lemma 5.1]{Cannizzaro20192DTriviality}, we have
\begin{equ}
  \label{eq:resolvent}
\tilde D(\lambda)\eqdef  \int_0^\infty e^{-\lambda t} \mathbf E[F_1(t)^2] \dd t= \frac2{\lambda^2}
\mathbb E[\phi(\omega) (\lambda-\mathcal G)^{-1}\phi(\omega) ],
\end{equ}
with $\phi$ defined in \eqref{eq:phi}. Hence, our analysis will focus
on the resolvent $(\lambda-\mathcal G)^{-1}$. Recall that $\mathbb E$ is the expectation with respect to the stationary law of the environment.


A first necessary step is to describe how $\mathcal G$ acts on elements in $L^2(\mathbb P)$, for which
we need a more accurate description of the latter space.
Since $\mathbb P$ is Gaussian (and given by the law of the curl of the smoothed Gaussian free field $\omega$),
$L^2(\mathbb P)$ admits a Wiener chaos decomposition which we now briefly
describe.

Let $H_0$ be the set containing constant random variables and $H_n$ be the closure of
the span of
\begin{equ}
  \label{eq:span}
  \left \{
  \psi=\sum_{j_1,\dots,j_n=1}^2\int_{\mathbb{R}^{2n}}f_\textbf{j}(x_{1:n}):\prod_{l=1}^n\omega_{j_l}(x_l):\dd x_{1:n}
  \right\},
\end{equ}
where $x_{1:n}$ is a short-hand notation for $(x_1,\dots,x_n)$, $\textbf{j}\eqdef(j_1,\dots,j_n)$,
$:\dots :$ denotes the Wick product associated to the measure $\mathbb P$
and the symmetric functions $f_\textbf{j}$'s are such that
\begin{equ}\label{eq:kernel}
  \hat{\psi}(p_{1:n})\eqdef\sum_{j_1,\dots,j_n=1}^2\prod_{k=1}^n\tilde{p}_{k,j_k}\hat{f}_\textbf{j}(p_{1:n})
\end{equ}
 satisfies
 \begin{equ}
   \label{modulo}
\int_{\mathbb{R}^{2n}}\prod_{i=1}^n\frac{\widehat{V}(p_i)}{|p_i|^2}|\hat \psi(p_{1:n})|^2\dd p_{1:n}<\infty\,.
\end{equ}
Above, $\hat f_\textbf{j}$ is the Fourier transform of $f_\textbf{j}$, $\tilde{p}_{k,1}=(p_{k})_2$ and
$\tilde{p}_{k,2}=-(p_{k})_1$, with $(p_k)_\ell$ the $\ell$-th component, $\ell=1,2$, of $p_k$.
Also, $\widehat V$ is the Fourier transform of the bump function $V$.
\begin{remark}
  The intuition behind \eqref{modulo} becomes clear upon noting that the components of
  $\omega$ are nothing but the derivatives of the smoothed Gaussian field $\xi$ (see~\eqref{eq:omegaxi}).
  Indeed, by performing an $n$-fold integration by parts in \eqref{eq:span},
  using the fact that, in Fourier space, $\partial_{(x_j)_\ell}$ corresponds to multiplication by $\iota (p_j)_\ell$
  (with $\iota=\sqrt{-1}$) and writing the covariance \eqref{eq:covxi} of the field $\xi$ in Fourier variables, one sees that \eqref{modulo} is just the $L^2(\mathbb P)$ norm squared of $\psi$ in \eqref{eq:span}.
\end{remark}

\begin{remark}
  \label{rem:phi}
  It is easy to see that the random variable
  $\phi(\omega)=\omega_1(0)$ belongs to $H_1$ and has kernel $\hat\phi(p)= p_2$.
\end{remark}

Then, by~\cite[Theorem 1.1.1]{Nualart}, $L^2(\mathbb P)$ can be orthogonally decomposed as
\begin{equ}
  \label{eq:Fock}
  L^2(\mathbb P)={\bigoplus_{n=0}^\infty H_n}
\end{equ}
and the expectation of the scalar product of $F,\,G\in L^2(\mathbb P)$ satisfies
\begin{equ}
\mathbb E [FG]=\sum_{n=1}^\infty \langle \psi_n,\phi_n\rangle\,.
\end{equ}
Above, $\psi_n$ and $\phi_n$ are the kernels of the projections of $F$ and $G$ onto $H_n$ and
the scalar product appearing at the right hand side is given by
\begin{equ}
  \label{eq:isometry}
  \langle \psi_n,\phi_n\rangle\eqdef n!\int_{\mathbb R^{2n}}\prod_{j=1}^n\frac{\widehat V(p_j)}{|p_j|^2}\,\overline{\hat \psi_n(p_{1:n})}\,\hat \phi_n(p_{1:n})\,\dd p_{1:n}\,.
\end{equ}

\begin{remark}
  \label{rem:abuseofnot}
  In what follows, we will implicitly identify a random variable $F$ in $H_n$ of the
  form \eqref{eq:span} with its kernel $\hat \psi_n$ in Fourier space,
  since this mapping is an isometry from $H_n$ to
 the set $L^2_{sym}(\mathbb R^{2n})$ of symmetric functions on $\mathbb R^{2n}$ endowed with the scalar product
 $\langle\cdot,\cdot\rangle$ in \eqref{eq:isometry}.
 In the same spirit, we will identify linear operators acting on $L^2(\mathbb P)$ with the corresponding
 linear operators acting on $\bigoplus_n L^2_{sym}(\mathbb R^{2n})$, and, with a slight abuse of notation,
 we will denote them using the same symbol.
\end{remark}

\begin{remark}
  \label{rem:differentconventions}
  With respect to~ \cite{Toth2012Superdiffusive2,Horvath2012Diffusive3} we are using different normalization conventions  in \eqref{eq:span} and in the scalar product in \eqref{eq:isometry}. More specifically, in the conventions of
  \cite{Toth2012Superdiffusive2,Horvath2012Diffusive3} there would be a factor $1/\sqrt{n!}$ in front of the integral in \eqref{eq:span} and no factor $n!$ in \eqref{eq:isometry}. In other words, our kernels $\psi_n$ equal those of \cite{Toth2012Superdiffusive2,Horvath2012Diffusive3} times $1/\sqrt{n!}$.
Our conventions are consistent with those of
\cite{Cannizzaro2020TheSuperdiffusivity} and of
\cite{Janson1997GaussianSpaces,Nualart};  we refer to these latter references for more details on Wiener chaos analysis.
\end{remark}

We are now ready to move back to the analysis of the generator
$\mathcal G$ of the environment process. As noted in~\cite{Toth2012Superdiffusive2},
$\mathcal G$ can be written as
\begin{equ}
\mathcal G = -\Delta+\genap-\genap^*,
\end{equ}
where $-\Delta$ and $\gena\eqdef\genap-\genap^*$ respectively denote the symmetric and
anti-symmetric part of $\mathcal G$ with respect to $\mathbb P$, and $\genap^*$ is
the adjoint of $\genap$ in $L^2(\mathbb P)$.
The action of $-\Delta$ and $\gena$ in Fock space is explicit.
First of all, $\Delta$ maps the $n$-th chaos $H_n$ into itself while $\genap$ (resp. $\genap^*$) maps
$H_n$ into $H_{n+1}$ (resp. $H_{n-1}$) and can therefore be interpreted as
a ``creation'' (resp. annihilation) operator. Moreover, $\Delta$ is diagonal in Fourier space
as it acts as a Fourier multiplier on the kernels,  while $\gena_+$ is not. Adopting the convention in
Remark \ref{rem:abuseofnot}, one has (see also~\cite[Section 2.1]{Toth2012Superdiffusive2})\footnote{If we had adopted the normalization conventions analogous to those of \cite{Horvath2012Diffusive3}, the factor $1/(n+1)$ in \eqref{eq:AplusF} would be replaced by $1/\sqrt{n+1}$. This is due to the different definition of the kernels, see Remark \ref{rem:differentconventions} above. }
\begin{equs}
\widehat{(-\Delta)\psi_n}(p_{1:n})&=|\sum_{i=1}^n p_i|^2\widehat\psi_n(p_{1:n})\\
\widehat{\genap\psi_n}(p_{1:n+1})&=\iota\frac{1}{n+1}\sum_{i=1}^{n+1}\Big(p_i\times\sum_{j=1}^{n+1}p_j\Big)\widehat{\psi_n}(p_{1:n+1\setminus i})\,, \label{eq:AplusF}
\end{equs}
$\iota=\sqrt{-1}$, for $\psi_n\in H_n$.
Above and throughout, we denote by $p_{1:n+1\setminus i}$ the collection $p_{1:n}=(p_1,\dots,p_{n+1})$ with $p_i$ removed.
Also, for $a,b$ two vectors in $\mathbb R^2$, by $a\times b$ we mean
the scalar given by the vertical component of the usual cross product $a\times b$, with $a,b$ viewed as vectors in $\mathbb R^3$. Explicitly,
$a\times b=|a||b|\sin\theta$, where $\theta$ is the angle between $a$ and $b$.

\begin{remark}
  Observe that the ``Laplacian'' $-\Delta$ is different from the one appearing in~\cite{Cannizzaro2020TheSuperdiffusivity},
  which acts instead as multiplication by $\sum_{i=1}^n |p_i|^2$.
  This represents a major technical difference which forces us to significantly modify the arguments therein.
  \end{remark}


At last, in light of the notations and conventions introduced above (see in particular Remark~\ref{rem:abuseofnot})
we rewrite \eqref{eq:resolvent} in Fock space as
\begin{equ}\label{eq:FockResolvent}
   \int_0^\infty e^{-\lambda t} \mathbf E[F_1(t)^2] \dd t= \frac2{\lambda^2}
\langle\phi, (\lambda-\mathcal G)^{-1}\phi\rangle,
\end{equ}
where $\phi$ is the random variable $\phi(\omega)$ in~\eqref{eq:phi} which lives in $H_1$ (see Remark~\ref{rem:phi}).

\section{The generator equation and the diffusivity}\label{sec:gen}

In order to obtain suitable bounds on the right hand side of \eqref{eq:FockResolvent} one should in principle
solve the generator equation $(\lambda-\mathcal G)\psi=\phi$ and then try to
evaluate $\langle \phi,\psi\rangle$.
While $\phi$ belongs to the first chaos, the operator $\mathcal G$ is not diagonal in the chaos decomposition
and finding $\psi$ explicitly is a rather challenging task.
A way out was first devised in \cite{landim2004superdiffusivity}. The idea is to truncate the generator $\mathcal G$
by defining $\mathcal G_n\eqdef I_{\le n}\mathcal G I_{\le n}$,
with $I_{\le n}$ the orthogonal projection onto $H_{\le n}\eqdef \oplus_{k\le n}H_k$ (the chaoses up to order $n$),
and then consider the solution $\psi^{(n)}\in H_{\le n}$ of the truncated generator equation
\begin{equ}
  \label{eq:troncata}
  (\lambda-\mathcal G_n)\,\psi^{(n)}=\phi.
\end{equ}
The advantage of this procedure is that it provides upper and lower
bounds (depending on the parity of $n$) on
\eqref{eq:FockResolvent}. Indeed, the following lemma,  which was first proven
in~\cite[Lemma 2.1]{landim2004superdiffusivity} (and whose proof straightforwardly carries out in the present case)
holds.

\begin{lem}\label{lemma:System}
  For every $n\ge 1$, one has
  \begin{equ}
    \label{eq:LQ}
    \langle\phi,\psi^{(2n)}\rangle \le \langle \phi, (\lambda-\mathcal G)^{-1}\phi\rangle=\langle\phi,\psi\rangle\le  \langle\phi, \psi^{(2n+1)}\rangle.
  \end{equ}
\end{lem}

The equation~\eqref{eq:troncata} coincides with the following hierarchical system of $n$ equations,
one for each component $\psi^{(n)}_k$ of $\psi^{(n)}$,
\begin{equ}\label{e:System}
\begin{cases}
\big(\lambda-\Delta\big)\psi^{(n)}_n-\genap\psi^{(n)}_{n-1}=0,\\
\big(\lambda-\Delta\big)\psi^{(n)}_{n-1}-\genap\psi^{(n)}_{n-2}+\genap^*\psi^{(n)}_{n}=0,\\
\dots\\
\big(\lambda-\Delta\big)\psi^{(n)}_1+\genap^*\psi^{(n)}_2=0\,.
\end{cases}
\end{equ}
Since $\phi$ belongs to the first chaos and different chaoses are orthogonal, in order to
estimate the terms at the left and right hand side of~\eqref{eq:LQ} we only
need to know $\psi^{(n)}_1$. The latter in turn can be obtained by solving the system~\eqref{e:System} iteratively
starting from $k=n$ so that we get
\begin{equ}
  \label{eq:LQ1}
  \langle \phi,\psi^{(n)}\rangle=  \langle \phi,\psi^{(n)}_1\rangle=  \langle\phi, (\lambda-\Delta+\Op_{n})^{-1}\phi\rangle
\end{equ}
where the self-adjoint operators $\Op_j$ are recursively defined as
\begin{equs}
  \label{eq:Opj}
  \Op_1&\eqdef 0\,,\\
  \Op_{j+1}&=\genap^*(\lambda-\Delta+\Op_j)^{-1}\genap\,,\qquad \text{for }j\geq 1\,.
\end{equs}
We remark that these operators are positive and leave each chaos invariant - that is
$\Op_j H_n\subset H_n$, for all $j,\,n\in\N$.


\subsection{Operator recursive estimates}

In view of \eqref{eq:LQ} and \eqref{eq:LQ1}, the proof of Theorem \ref{th:main} must entail
a good understanding of the operators $\Op_j$'s in~\eqref{eq:Opj}.
In particular, we need to derive suitable (upper and lower) bounds on them and this is the
content of the main result of this section, Theorem~\ref{th:bounds}.
To state it, we need a few preliminary definitions.

For $k\in\mathbb{N}$, $x>0$ and $z\ge0$ we define $\Ll$, ${\LB}_k$ and $\UB_k$ as follows
\begin{equs}
    \Ll(x,z)&= z+\log(1+x^{-1}),\label{eq:defL}\\
    {\LB}_k(x,z)&=\sum_{0\le j\leq k}\frac{(\frac{1}{2}\log \Ll(x,z))^j}{j!} \quad\text{ and }\quad {\UB}_k(x,z)=\frac{\Ll(x,z)}{{\LB}_k(x,z)}\label{eq:defLBUB}
\end{equs}
and for $k\geq 1$, $\sigma_k$, as
\begin{equ}
   \sigma_k(x,z) =
   \begin{cases}
\UB_\frac{k-2}{2}(x,z),\quad\text{ if $k$ is even,}\\
\LB_\frac{k-1}{2}(x,z),\quad \text{ if $k$ is odd.}\\
\end{cases}
\end{equ}
Note that $\sigma_1\equiv 1$. All the  properties we need on the functions $\UB_k,\LB_k$ are summarized in Lemma \ref{lemma:LBUBproperties}.
Further let
\begin{equ}
  \label{eq:K1K2}
    z_k(n)=K_1(n+k)^{2+2\eps} \quad\text{ and } \quad f_k(n)=K_2\sqrt{z_k(n)},
\end{equ}
where $K_1$, $K_2$ are absolute constants (chosen sufficiently large, 
so that~\eqref{e:CondK1K21},~\eqref{e:CondK1K22} and~\eqref{e:CondK1K23} below, hold) and
$\eps$ is the small positive constant that appears in the statement of Theorem~\ref{th:main}.

Finally, for $k\geq1$ let $\mathcal{S}_k$ be the operator whose multiplier is $\sigma_k$, i.e.
\begin{equ}
    \mathcal{S}_k=
    \begin{cases}
    f_k(\mathcal{N})\,\sigma_k(\lambda-\Delta,z_k(\mathcal{N})) &\text{ if $k$ is even,}\\
    \frac{1}{f_k(\mathcal{N})}\big(\sigma_k(\lambda-\Delta,z_k(\mathcal{N}))-f_k(\mathcal{N})\big) &\text{ if $k$ is odd,}\\
    \end{cases}
\end{equ}
where $\mathcal{N}$ is the number operator acting on the $n$-th chaos as multiplication by $n$,
i.e. $(\mathcal{N}\phi_n)=n\phi_n$ for $\phi_n\in H_n$.
We are now ready to state the following theorem.

\begin{thm}
  \label{th:bounds}
  For any $\eps>0$, the constants $K_1,K_2$ in \eqref{eq:K1K2} can be chosen in such a way that the following holds.
  For $0<\lambda\le 1$ and $k\ge1$, one has the operator bounds
 \begin{equ}
    \label{eq:LBOp}
    \Op_{2k-1}\ge c_{2k-1}\,(-\Delta)\cS_{2k-1}
  \end{equ}
  and
  \begin{equ}
    \label{eq:UBOp}
    \Op_{2k}\le c_{2k}\,(-\Delta)\cS_{2k}
  \end{equ}
  where $c_{1}=1$ and
  \begin{equ}
    \label{eq:c}
    c_{2k}=\frac{\pi}{c_{2k-1}}\left(1+\frac1{k^{1+\eps}}\right),\qquad
    c_{2k+1}=\frac{\pi}{c_{2k}}\left(1-\frac1{(k+1)^{1+\eps}}\right).
  \end{equ}
\end{thm}

\begin{remark}\label{rem:Limitcj}
A crucial aspect we need to stress is that, as $j\to\infty$, $c_{2j}$  tends to a {\it finite} constants larger than $1$, while  $c_{2j+1}$ tends to a {\it strictly positive} constant smaller than $1$.
\end{remark}

\subsection{Generalities about the operators}

In this section we collect some preliminary facts and bounds concerning
operators in Fock space. In all the statements herein, $S$ will be a {\it diagonal operator},
meaning that $S$ commutes with $\mathcal N$ (that is, it maps the
$n$-th chaos $H_n$ into itself) and is diagonal in the Fourier
basis, i.e. it acts in Fourier space as multiplication by a
function of the momenta. The Fourier multiplier of $S$ will be denoted $\mathfrak s$,
and actually $\mathfrak s$ is the collection $(\mathfrak s_n)_{n\ge1}$, with $\mathfrak s_n$ the Fourier multiplier on $H_n$.
It is understood that $\mathfrak s_n$ is a symmetric function of its $n$ arguments.

\begin{lem}
  \label{lemma:splitting}
  Let $S$ be a positive diagonal operator and let $\mathfrak s$ be its Fourier multiplier. For any element $\psi$ of  $H_n$ we can write
  \begin{equ}
    \langle\psi, \genap^* S \genap \psi\rangle=     \langle\psi, \genap^* S \genap \psi\rangle_{\mathrm{Diag}}+
        \langle\psi, \genap^* S \genap \psi\rangle_{\mathrm{ Off}}
  \end{equ}
  where the ``diagonal part'' is defined as
  \begin{equs}[eq:diag]
    \langle\psi, &\genap^* S \genap \psi\rangle_{\mathrm{Diag}}\\
    &\eqdef n!\int_{\mathbb R^{2(n+1)}}\prod_{j=1}^{n+1} \frac{\widehat V(p_j)}{|p_j|^2}|\hat\psi(p_{1:n})|^2\mathfrak s_{n+1}(p_{1:n+1})\left(
    p_{n+1}\times \sum_{j=1}^np_j\right)^2\dd p_{1:n+1}
\end{equs}
while the ``off-diagonal part'' is
\begin{equs}[eq:offdiag]
  \langle\psi, &\genap^* S \genap \psi\rangle_{\mathrm{ Off}}\\
  &\eqdef n!n \int_{\mathbb R^{2(n+1)}}\prod_{j=1}^{n+1} \frac{\widehat V(p_j)}{|p_j|^2}\overline{\hat \psi(p_{1:n})}\hat\psi(p_{1:n+1\setminus n}) \mathfrak s_{n+1}(p_{1:n+1})\times\\
&\qquad\qquad\qquad\qquad\qquad\times\left(p_{n+1}\times \sum_{i=1}^{n+1}p_i\right) \left(p_n\times  \sum_{i=1}^{n+1}p_i\right)
 \dd p_{1:n+1}.
\end{equs}
\end{lem}
\begin{proof}
Expanding the inner product using (\ref{eq:AplusF}) we obtain:
\begin{align*}
&\langle \genap\psi_n, S\genap\psi_n\rangle=\\
&\frac{(n+1)!}{(n+1)^2}\int_{\mathbb{R}^{2(n+1)}}\prod_{i=1}^{n+1}\frac{\widehat{V}(p_i)}{|p_i|^2}
\mathfrak s_{n+1}(p_{1:n+1})\left|\sum_{i=1}^{n+1}\hat{\psi}(p_{1:n+1\setminus i})\Big(p_i\times\sum_{j=1}^{n+1}p_j\Big)\right|^2 \dd p_{1:n+1}.
\end{align*}
The ``diagonal'' and ``off-diagonal'' refer to the squared sum.
The former is the contribution of the squared summands while the latter comes from all the cross terms.
Hence, the diagonal part is
\begin{align*}
\frac{n!}{n+1}&\int_{\mathbb{R}^{2(n+1)}}\prod_{i=1}^{n+1}\frac{\widehat{V}(p_i)}{|p_i|^2}
\mathfrak s_{n+1}(p_{1:n+1})\sum_{i=1}^{n+1}|\hat{\psi}(p_{1:n+1\setminus i})|^2\Big(p_i\times\sum_{j=1}^{n+1}p_j\Big)^2 \dd p_{1:n+1}\\
&=n!\int_{\mathbb R^{2(n+1)}}\prod_{j=1}^{n+1} \frac{\widehat V(p_j)}{|p_j|^2}|\hat\psi(p_{1:n})|^2\mathfrak s_{n+1}(p_{1:n+1})\Big(
p_{n+1}\times \sum_{j=1}^np_j\Big)^2\dd p_{1:n+1},
\end{align*}
where we pulled out the sum and used that $\hat\psi$ is symmetric in its arguments.
For the off-diagonal part, one follows the same procedure.
Since there are in total $n(n+1)$ summands, a factor $n$ is left in front of the integral.
\end{proof}

The next two results will be used in the bounds on the diagonal and off-diagonal parts respectively.
In order to appreciate them, note that at the right hand side of both~\eqref{eq:diag}
and~\eqref{eq:offdiag}, there appears the vector product.
\begin{lem}
  \label{lemma:sigmasigma'}
  Let $S$ be a positive diagonal operator, and let $\mathfrak s$ be the associated Fourier multiplier.
  If for every integer $n$ and for every $p_{1:n}\in \mathbb R^{2n}$ with $\sum_{k=1}^n p_k\ne 0$
  \begin{equ}[e:se]
    \int_{\mathbb R^2}\widehat V(q)(\sin\theta)^2\mathfrak s_{n+1}(p_{1:n},q)\dd q\le \tilde{\mathfrak s}_n(p_{1:n})
  \end{equ}
  with $\theta$ the angle between $q$ and $\sum_{k=1}^n p_k$, then for every  $\psi$
  \begin{equ}[e:allora]
    \langle \psi, \genap^* S \genap \psi\rangle_{\mathrm{Diag}}\le \langle \psi,(-\Delta) \tilde S\psi\rangle
  \end{equ}
  where $\tilde S$ is the diagonal operator whose Fourier multiplier is $\tilde{\mathfrak s}$.

If the inequality in \eqref{e:se} is reversed, then \eqref{e:allora} holds with the reversed inequality.
\end{lem}

\begin{proof}
Starting from (\ref{eq:diag}) and denoting $q=p_{n+1}$ we get that the left-hand side equals:
\begin{equs}
n!\int_{\mathbb R^{2n}}&\prod_{j=1}^{n} \frac{\widehat V(p_j)}{|p_j|^2}|\hat\psi(p_{1:n})|^2|\sum_{k=1}^np_k|^2\int_{\mathbb R^2}\mathfrak s_{n+1}(p_{1:n},q)\widehat{V}(q)(\sin\theta)^2\dd q\dd p_{1:n}\\
&\leq n!\int_{\mathbb R^{2n}}\prod_{j=1}^{n} \frac{\widehat V(p_j)}{|p_j|^2}|\hat\psi(p_{1:n})|^2|\sum_{k=1}^np_k|^2\tilde{\mathfrak s}_n(p_{1:n})=\langle \psi,(-\Delta) \tilde S\psi\rangle\,,
\end{equs}
where we used that $a\times b=|a||b|\sin\theta$, $\theta$ being the angle between $a$ and $b$, and~\eqref{e:se}.
Since every step except the assumption is an equality, the other direction also holds.
\end{proof}

\begin{lem}
  \label{lemma:off} Let $S$ be a diagonal, positive operator with Fourier multiplier $\mathfrak s$. If for every integer $n$ and every $p_{1:n}\in\mathbb R^{2n} $ one has
  \begin{equ}
    |\sum_{i=1}^np_i|\int_{\mathbb R^2}\widehat V(q)\frac{(\sin\theta)^2\mathfrak s_{n+1}(p_{1:n},q)}{|q+\sum_{i=1}^{n-1}p_i|}\dd q\le \tilde {\mathfrak s}_n(p_{1:n})
  \end{equ}
  with $\theta$ the angle between $q$ and $\sum_{i=1}^n p_i$, then for
  every $n\in\mathbb N,\psi\in H_n$ one has
  \begin{equ}
   |  \langle\psi, \genap^* S \genap \psi\rangle_{\mathrm{ Off}}|\le n \langle \psi, (-\Delta)\tilde S\psi\rangle,
  \end{equ}
  with $\tilde S$ the diagonal operator of Fourier multiplier $\tilde{\mathfrak s}$.
\end{lem}
\begin{proof}
We start by bounding the left-hand side of (\ref{eq:offdiag}) as
\begin{equs}
n!n \int_{\mathbb R^{2(n+1)}}&\prod_{j=1}^{n+1} \frac{\widehat V(p_j)}{|p_j|^2}|
\hat \psi(p_{1:n})||\hat\psi(p_{1:n+1\setminus n})|
\mathfrak s_{n+1}(p_{1:n+1})\\
&\qquad\qquad\qquad\times
\Big|p_{n+1}\times \sum_{k=1}^np_k\Big| \Big|p_n\times  \Big(\sum_{k=1}^{n-1}p_k+p_{n+1}\Big)\Big|
\dd p_{1:n+1}\\\label{tobebounded}
&=
n!n \int_{\mathbb R^{2(n+1)}}\prod_{j=1}^{n+1} \frac{\widehat V(p_j)}{|p_j|^2}\mathfrak s_{n+1}(p_{1:n+1})
\Phi(p_{1:{n+1}})\Phi(p_{1:{n-1}},p_{n+1},p_n)\\
&\qquad\qquad\qquad\qquad\qquad\qquad\qquad\times
|\sum_{k=1}^{n-1}p_k+p_{n+1}||\sum_{k=1}^np_k|
\dd p_{1:n+1}
\end{equs}
where \[\Phi(p_{1:{n+1}})=\frac{|\hat \psi(p_{1:n})||p_{n+1}\times \sum_{k=1}^np_k|}{|\sum_{k=1}^{n-1}p_k+p_{n+1}|}.\]
We apply Cauchy-Schwartz and exploit symmetry of $\psi$  to bound \eqref{tobebounded} from above by
\begin{equ}
n!n \int_{\mathbb R^{2(n+1)}}\prod_{j=1}^{n+1} \frac{\widehat V(p_j)}{|p_j|^2}\mathfrak s_{n+1}(p_{1:n+1})
\Phi(p_{1:{n+1}})^2
|\sum_{k=1}^{n-1}p_k+p_{n+1}||\sum_{k=1}^np_k|
\dd p_{1:n+1}.
\end{equ}
Now set $s_1=\sum_{k=1}^np_k$, $s_2=\sum_{k=1}^{n-1}p_k$ and $q=p_{n+1}$, which gives
\begin{equs}
n!n &\int_{\mathbb R^{2(n+1)}}\prod_{j=1}^{n+1} \frac{\widehat V(p_j)}{|p_j|^2}\mathfrak s_{n+1}(p_{1:n+1})
\frac{|\hat \psi(p_{1:n})|^2|q\times s_1|^2}{|q+s_2|}
|s_1|
\dd p_{1:n+1}\\
&= n!n\int_{\mathbb R^{2n}}\prod_{j=1}^{n}\left( \frac{\widehat V(p_j)}{|p_j|^2}\dd p_j\right)|\hat \psi(p_{1:n})|^2|s_1|^3
\int_{\mathbb R^2}\mathfrak s_{n+1}(p_{1:n},q)
\frac{\widehat{V}(q)(\sin\theta)^2}{|q+s_2|}\dd q
\\
&\leq
n!n\int_{\mathbb R^{2n}}\prod_{j=1}^{n} \frac{\widehat V(p_j)}{|p_j|^2}|\hat \psi(p_{1:n})|^2|s_1|^2\tilde{\mathfrak s}_n(p_{1:n})\dd p_{1:n}
=n \langle \psi, (-\Delta)\tilde S\psi\rangle\,,
\end{equs}
which concludes the proof.
\end{proof}

\subsection{Proof of Theorem \ref{th:bounds}}

This section is devoted to Theorem \ref{th:bounds}.
We will first show the lower bound and then the upper bound, both by induction on $k$. The induction switches from lower to upper bounds and viceversa, as follows:
For $k=1$ the bound \eqref{eq:LBOp} will be trivial; given  \eqref{eq:LBOp} for $k=1$ we will deduce \eqref{eq:UBOp} with $k=1$, then \eqref{eq:LBOp} with $k=2$ and so on.

\begin{proof}[Proof of the lower bound~\eqref{eq:LBOp}]
 For $k=1$ ~\eqref{eq:LBOp} trivially holds as $\mathcal H_1$  is by definition zero while $\mathcal S_1$ is non-positive if the constant $K_2$ in the definition \eqref{eq:K1K2} is large enough.

We need then to prove \eqref{eq:LBOp} with $k\ge1$ and $2k-1$ replaced by $2k+1$. Assume by induction that~\eqref{eq:UBOp} holds. Then, we have
\begin{equ}
  \label{eq:osserviamo}
    \mathcal{H}_{2k+1}=\genap^*(\lambda-\Delta+\mathcal{H}_{2k})^{-1}\genap\geq \genap^*(\lambda-\Delta(1+c_{2k}\cS_{2k}))^{-1}\genap\,.
  \end{equ}
  For $\psi\in H_n$, we apply Lemma \ref{lemma:splitting} with
  $S=(\lambda-\Delta(1+c_{2k}\cS_{2k}))^{-1}$ and we split
  \begin{equ}
    \label{eq:split}
\langle \psi,\genap^*(\lambda-\Delta(1+c_{2k}\cS_{2k}))^{-1}\genap \psi\rangle
\end{equ}
into diagonal and off-diagonal part. In order to control the
former from below, we exploit Lemma~\ref{lemma:sigmasigma'} according to which it suffices to consider
\begin{equ}
  \label{eq:numden}
    \int_{\mathbb{R}^2}\frac{\widehat{V}(q)(\sin\theta)^2}{\lambda+|p+q|^2(1+c_{2k}f_{2k}(n+1)\UB_{k-1}(\lambda+|p+q|^2,z_{2k}(n+1)))}\dd q,
\end{equ}
where $p=\sum_{i=1}^n p_i\ne0$ and $\theta $ is the angle between $p$ and $q$.
Note that the functions $f_{2k},z_{2k}$ have argument $n+1$ because $\genap \psi\in H_{n+1}$
but, by~\eqref{eq:K1K2}, $f_{2k}(n+1)=f_{2k+1}(n), z_{2k}(n+1)=z_{2k+1}(n)$.
To lighten the notation, throughout the proof we will omit the argument $n$ and
write $z_{2k+1},f_{2k+1}$ instead of $z_{2k+1}(n),f_{2k+1}(n)$.

The denominator in \eqref{eq:numden} is upper bounded by
\begin{equ}
    c_{2k}f_{2k+1}\left(1+\frac{1}{f_{2k+1}}\right)(\lambda+|p+q|^2\UB_{k-1}(\lambda+|p+q|^2,z_{2k+1})),
\end{equ}
as $c_{2k},f_{2k+1}$ and $\UB_{k-1}$ are all larger than one. Thus we can concentrate on
\begin{equ}
    \int_{\mathbb{R}^2}\frac{\widehat{V}(q)(\sin\theta)^2}{\lambda+|p+q|^2\UB_{k-1}(\lambda+|p+q|^2,z_{2k+1})}\dd q.
\end{equ}
For this we first apply Lemmas \ref{lemma:replacement} and~\ref{lemma:ff'} to obtain the lower bound
\begin{equs}
    &\frac\pi2\int_{\lambda+|p|^2}^1\frac{d\rho}{\rho \,\UB_{k-1}(\rho,z_{2k+1})}-
    C_{\Di}\frac{\LB_{k-1}(\lambda+|p|^2,z_{2k+1})}{\sqrt{z_{2k+1}}}\\
    &\geq\frac\pi2\int_{\lambda+|p|^2}^1\frac{d\rho}{(\rho+\rho^2) \UB_{k-1}(\rho,z_{2k+1})}-
    C_{\Di}\frac{\LB_{k}(\lambda+|p|^2,z_{2k+1})}{\sqrt{z_{2k+1}}}\,.
\end{equs}
From Lemma \ref{lemma:LBUBproperties} we have that the primitive of the integrand is $-2\LB_k(\rho,z_{2k+1})$, so that the last expression equals
\begin{equs}
    \pi \LB_{k}&(\lambda+|p|^2,z_{2k+1})-\pi \LB_{k}(1,z_{2k+1})-
    C_{\Di}\frac{\LB_{k}(\lambda+|p|^2,z_{2k+1})}{\sqrt{z_{2k+1}}}
\\
    &\geq \pi \LB_{k}(\lambda+|p|^2,z_{2k+1})-\pi \frac{f_{2k+1}}2-
    C_{\Di}\frac{\LB_{k}(\lambda+|p|^2,z_{2k+1})}{\sqrt{z_{2k+1}}}\,,
\end{equs}
where in the first inequality we need to choose $K_2$ large enough in
\eqref{eq:K1K2} so that for all $k$ and $n$,
\begin{equ}[e:CondK1K21]
    \LB_k(1,z_{2k+1})\leq\sqrt{\Ll(1,z_{2k+1})}=\sqrt{\log(2)+z_{2k+1}}\leq \frac12 f_{2k+1}
  \end{equ}
  (see also \eqref{e:bdonL}).
Altogether, the diagonal part of \eqref{eq:split} is lower bounded as
$\langle \psi, (-\Delta)\tilde S\psi\rangle$, with
\begin{equ}[eq:alto]
\tilde S=\left(1+\frac1{f_{2k+1}(1)}\right)^{-1}  \frac{\pi}{c_{2k}}\left[\frac{\LB_k(\lambda-\Delta,z_{2k+1}(\mathcal N))}{f_{2k+1}(\mathcal N)}\left(1-\frac{C_{\Di}}{\pi\sqrt{z_{2k+1}(1)}}\right) -\frac12
    \right]
\end{equ}
(in two instances, we have lower bounded $z_{2k+1}=z_{2k+1}(n)$ and $f_{2k+1}=f_{2k+1}(n)$ with the same quantities for $n=1$).

For the off-diagonal terms in \eqref{eq:split} we use Lemma \ref{lemma:off} so that, calling $p:=\sum_{i=1}^n p_i$ and $p':=\sum_{i=1}^{n-1}p_i$, we have to upper bound
\begin{equ}[e:OffDiagLB]
n |p|\int_{\mathbb R^2}\frac{\widehat V(q)(\sin\theta)^2}{(\lambda+|p+q|^2(1+c_{2k}f_{2k+1}\UB_{k-1}(\lambda+|p+q|^2,z_{2k+1})))|p'+q|}\dd q\,.
\end{equ}
Thanks to Lemmas \ref{lemma:off-diagonals} and \ref{lemma:ff'}, applied with $f(x,z)=c_{2k}f_{2k+1}\UB_{k-1}(x,z)$ and $g(x,z)=\frac{1}{c_{2k}f_{2k+1}}\LB_{k-1}(x,z)$,
this expression is upper bounded by
\begin{equs}
    & \frac{n\, C_{\oDi}}{c_{2k}f_{2k+1}z_{2k+1}} \LB_{k-1}(\lambda+|p|^2,z_{2k+1})\notag\\
  &\qquad\leq \frac{n\,C_{\oDi}}{c_{2k}f_{2k+1} z_{2k+1}}\LB_{k}(\lambda+|p|^2,z_{2k+1})\notag\\
  \label{eq:alto'}
    &\qquad\leq \frac{C_{\oDi}}{c_{2k}f_{2k+1} }\frac1{K_1(2k+1)^{1+\eps}}\LB_{k}(\lambda+|p|^2,z_{2k+1}),
\end{equs}
where we used monotonicity properties of $\LB_k$, the definition of $z_{2k+1}=z_{2k+1}(n)$ in~\eqref{eq:K1K2}
and in particular the fact that \[\frac{n}{z_{2k+1}(n)}= \frac{n}{K_1(2k+1+n)^{2+2\eps}}\leq \frac{1}{K_1(2k+1+n)^{1+\eps}}.\]

Combining \eqref{eq:alto} and \eqref{eq:alto'}, together with Lemmas \ref{lemma:sigmasigma'} and \ref{lemma:off}, we conclude that $\genap^* (\lambda-\Delta(1+c_{2k}\cS_{2k}))^{-1}\genap$ is lower bounded by
\begin{equ}
  (-\Delta)\frac{\pi}{c_{2k}}\left[\frac{\LB_k(\lambda-\Delta,z_{2k+1}(\mathcal N))}{f_{2k+1}(\mathcal N)}A-B\right]
\end{equ}
where $A$ and $B$ are given by
\begin{equs}
A&=\left(1-\frac{C_{\Di}}{\pi\sqrt{z_{2k+1}(1)}}\right)\left(1+\frac1{f_{2k+1}(1)}\right)^{-1}-\frac{C_{\oDi}}{\pi K_1(2k+1)^{1+\eps}}
\\
B&=\frac12\left(1+\frac1{f_{2k+1}(1)}\right)^{-1}
\end{equs}
and thanks to \eqref{eq:osserviamo} the same lower bound holds for $\Op_{2k+1}$.
Note that, provided the constants $K_1,K_2$ in \eqref{eq:K1K2} are large, one has
\begin{equ}[e:CondK1K22]
A\ge  1-\frac1{(k+1)^{1+\eps}}, \qquad B\le  1-\frac1{(k+1)^{1+\eps}}\,.
\end{equ}
Therefore, we have proven \eqref{eq:LBOp} (with $2k+1$ instead of $2k-1$) with $c_{2k+1}$ given by \eqref{eq:c}.
\end{proof}

\begin{proof}[Proof of the upper bound~\eqref{eq:UBOp}]
 For $k\ge 1$,
again by the induction hypothesis we have
\begin{equ}
    \mathcal{H}_{2k}=\genap^*(\lambda-\Delta+\mathcal{H}_{2k-1})^{-1}\genap\leq \genap^*(\lambda-\Delta(1+c_{2k-1}\cS_{2k-1}))^{-1}\genap\,.
  \end{equ}
We split $\langle\psi,\genap^*(\lambda-\Delta(1+c_{2k-1}\cS_{2k-1}))^{-1}\genap\psi\rangle$ into diagonal and off-diagonal parts as in Lemma \ref{lemma:splitting}.
By Lemma \ref{lemma:sigmasigma'} for the diagonal part we need to upper-bound the integral
\begin{equs}
    \int_{\mathbb{R}^2}&\frac{\widehat{V}(q)(\sin\theta)^2}{\lambda+|p+q|^2(1+\frac{c_{2k-1}}{f_{2k}}(\LB_{k-1}(\lambda+|p+q|^2,z_{2k})-f_{2k}))}\dd q\\
    &\qquad\qquad\leq \frac{f_{2k}}{c_{2k-1}}\int_{\mathbb{R}^2}\frac{\widehat{V}(q)(\sin\theta)^2}{\lambda+|p+q|^2\LB_{k-1}(\lambda+|p+q|^2,z_{2k})}\dd q
\end{equs}
where we used $f_{2k-1}(n+1)=f_{2k}(n)$, the same for $z$ (and we suppressed the argument of both)
and, in the second step, exploited the fact that $c_{2k-1}<1$ and $f_{2k}>1$.
By Lemmas \ref{lemma:replacement} and \ref{lemma:ff'}, the latter is bounded above by
\begin{align*}
\frac{f_{2k}\pi}{2c_{2k-1}}\left(\int_{\lambda+|p|^2}^1\frac{\dd\rho}{\rho \LB_{k-1}(\rho,z_{2k})}+\frac{C_\text{diag}\UB_{k-1}(\lambda+|p|^2,z_{2k})}{\sqrt{z_{2k}}}\right)\,.
\end{align*}
The integral can be controlled via Lemma \ref{lemma:rhosquared}, so that
\begin{align*}
\int_{\lambda+|p|^2}^1\frac{\dd\rho}{\rho \LB_{k-1}(\rho,z_{2k})}&\leq\int_{\lambda+|p|^2}^1\frac{\dd\rho}{(\rho +\rho^2)\LB_{k-1}(\rho,z_{2k})}
+C\frac{\UB_{k-1}(\lambda+|p|^2,z_{2k})}{z_{2k}}\\
&\leq2\UB_{k-1}(\lambda+|p|^2,z_{2k})+C\frac{\UB_{k-1}(\lambda+|p|^2,z_{2k})}{z_{2k}}\,,
\end{align*}
the last passage being a consequence of Lemma~\ref{lemma:LBUBproperties}.

For the off-diagonal terms, we argue as in the analysis of~\eqref{e:OffDiagLB}, so that we need to control
\begin{equ}
n |p|\int_{\mathbb R^2}\frac{\widehat V(q)(\sin\theta)^2}{(\lambda+|p+q|^2(1+\frac{c_{2k-1}}{f_{2k}}(\LB_{k-1}(\lambda+|p+q|^2,z_{2k})-f_{2k})))|p'+q|}\dd q\,.
\end{equ}
Once again, we can pull out the factor $\frac{f_{2k}}{c_{2k-1}}$ and apply once more Lemmas~\ref{lemma:off-diagonals} and \ref{lemma:ff'}, this time with $f(x,z)=\LB_{k-1}(x,z)$ and $g(x,z)=\UB_{k-1}(x,z)$. Hence we obtain
\begin{equs}
\frac{f_{2k}}{c_{2k-1}}\frac{nC_{\oDi}\UB_{k-1}(\lambda+|p|^2),z_{2k})}{z_{2k}}\leq \frac{f_{2k}}{c_{2k-1}}\frac{C_{\oDi}\UB_{k-1}(\lambda+|p|^2,z_{2k})}{K_1(n+2k)^{1+2\eps}}.
\end{equs}
%
%
%
Collecting these upper bounds and using the fact that $z_{2k}(n)>z_{2k}(1)$,
we conclude that $\genap^* (\lambda-\Delta(1+c_{2k-1}\cS_{2k-1}))^{-1}\genap$ is upper bounded by
\begin{equ}
\frac{\pi}{c_{2k-1}}A' (-\Delta)\cS_{2k}
\end{equ}
where this time, upon choosing $K_1$ big enough, we have
\begin{equ}[e:CondK1K23]
A'=1+\frac{C_{\Di}}{\pi\sqrt{K_1}(2k)^{1+\eps}}+
\frac{C}{\pi K_1(2k)^{2+2\eps}}+
\frac{C_{\oDi}}{\pi K_1(2k)^{1+2\eps}}\leq 1+\frac{1}{k^{1+\eps}}\,.
\end{equ}
It follows that~\eqref{eq:UBOp} holds with $c_{2k}$ satisfying~\eqref{eq:c}.

Let us remark that constants $K_1$ and $K_2$ such that~\eqref{e:CondK1K21},~\eqref{e:CondK1K22} and~\eqref{e:CondK1K23} hold for all $k,\,n\in\N$ clearly exist, so that the proof of Theorem~\ref{th:bounds} is concluded.
\end{proof}

\section{Proof of Theorem \ref{th:main}}\label{sec:Proof}

This section is devoted to the proof of Theorem~\ref{th:main} and
shows how to exploit the iterative bounds derived in the previous
section. Recall from Section \ref{sec:prelim} that it suffices to prove \eqref{eq:main} with $D(\lambda)$ replaced by $\tilde D(\lambda)$ defined in \eqref{eq:resolvent}.

\begin{proof}[Proof of Theorem \ref{th:main}]
Let us begin with the upper bound.
By Lemma \ref{lemma:System} and~\eqref{eq:FockResolvent}, we have
\begin{equ}
\frac{\lambda^2}2 \tilde D(\lambda)\leq \langle\phi,\psi^{(2k+1)}\rangle=\langle\phi,(\lambda-\Delta+\mathcal H_{2k+1})^{-1}\phi\rangle,
\end{equ}
for $\phi$ such that $\hat\phi(q)=q_2 $ (see Remark~\ref{rem:phi}),
which in turn, by Theorem \ref{th:bounds}, is bounded above by
\begin{equs}
\langle\phi,&(\lambda-\Delta(1+c_{2k+1}\mathcal S_{2k+1}))^{-1}\phi\rangle\\
&=
\int_{\mathbb R^2}
\frac{\widehat{V}(q)}{|q|^2}
\frac{|\hat\phi(q)|^2\dd q}{\lambda+|q|^2(1+\frac{c_{2k+1}}{f_{2k+1}}(\LB_k(\lambda+|q|^2,z_{2k+1})-f_{2k+1}))}\\
&\leq
\frac{f_{2k+1}}{c_{2k+1}}\int_{\mathbb R^2}
\widehat{V}(q)
\frac{\dd q}{\lambda+|q|^2\LB_k(\lambda+|q|^2,z_{2k+1})}\,.\label{e:almostthere}
\end{equs}
Note that, as $\phi\in H_1$, $f_{2k+1}$ and $z_{2k+1}$ are $f_{2k+1}(1)$ and $z_{2k+1}(1)$, that is, are constants depending only on $k$.

In view of~\eqref{eq:c}, we can replace $c_{2k+1}$ with its $k\to\infty$ limit.
By Eq. \eqref{e:ApproxDobbiamoFinire} in Lemma \ref{lemma:replacement},~\eqref{e:almostthere} is controlled,
up to a multiplicative absolute constant, by
\begin{equs}[e:UBcomplete]
 f_{2k+1}&\left[\int_\lambda^1\frac{\dd \rho}{\rho\LB_k(\rho,z_{2k+1})}+\frac{\UB_k(\lambda,z_{2k+1})}{\sqrt{z_{2k+1}}}\right]
   \\
   & \lesssim f_{2k+1}\left[\int_\lambda^1\frac{\dd \rho}{(\rho+\rho^2)\LB_k(\rho,z_{2k+1})}+\frac{\UB_k(\lambda,z_{2k+1})}{\sqrt{z_{2k+1}}}\right]\\
&\lesssim f_{2k+1}\UB_k(\lambda,z_{2k+1})\lesssim
f_{2k+1}\frac{\Ll(\lambda,0)+z_{2k+1}}{\LB_k(\lambda,0)},
\end{equs}
where in the first inequality we applied Lemma~\ref{lemma:rhosquared},
in the second Lemma~\ref{lemma:LBUBproperties} and in the last the monotonicity
of $\LB_k(\cdot,z)$ with respect to $z$.
We now recall that the central limit theorem, applied to Poisson random variables of rate one, gives that
\begin{equ}[e:CLT]
\lim_{k\to\infty}\sum_{j=0}^{k}\frac{k^j}{j!} e^{-k}=\frac12\,.
\end{equ}
Hence, by choosing
\begin{equ}[e:k]
k=k(\lambda)=\Big\lfloor\frac{\log \Ll(\lambda,0)}{2}\Big\rfloor.
\end{equ}
in~\eqref{e:CLT} and recalling the definition of $\LB_k$ in~\eqref{eq:defLBUB},
we have that for $\lambda$ sufficiently small
\begin{equ}[e:UBFinal]
\frac{e^{-k}}{\LB_k(\lambda,0)e^{-k}}\lesssim \frac{1}{\sqrt{\Ll(\lambda,0)}}\,.
\end{equ}
Plugging this  into~\eqref{e:UBcomplete} and using the definition of $z_{2k+1}=z_{2k+1}(1)$ and
$f_{2k+1}$ in \eqref{eq:K1K2},
we ultimately get the upper bound
\begin{equ}
\lambda^2\tilde D(\lambda)\lesssim (\log\Ll(\lambda,0))^{1+\eps}\sqrt{\Ll(\lambda,0)}
\end{equ}
which is the desired one, since \[\Ll(\lambda,0)=\log\left(1+\frac1\lambda\right)\stackrel{\lambda\to0}\sim |\log \lambda|.\]

For the lower bound, we argue similarly. Again by Lemma \ref{lemma:System}, we have
\begin{equ}
\frac{\lambda^2}{2}\tilde D(\lambda)\geq \langle\phi,\psi^{(2k)}\rangle=\langle\phi,(\lambda-\Delta+\mathcal H_{2k})^{-1}\phi\rangle,
\end{equ}
for $\phi$ such that $\hat\phi(q)=q_2 $,
which in turn, by Theorem \ref{th:bounds}, is bounded below by
\begin{equs}
  \langle\phi,(\lambda-\Delta(1+c_{2k}\mathcal S_{2k}))^{-1}\phi\rangle&\geq\int_{\mathbb R^2}
\frac{\widehat{V}(q)}{|q|^2}
\frac{|\hat\phi(q)|^2\dd q}{\lambda+|q|^2(1+c_{2k}f_{2k}\UB_{k-1}(\lambda+|q|^2,z_{2k}))}\\
&\gtrsim\frac{1}{f_{2k}}\int_{\mathbb R^2}
\frac{\widehat{V}(q)}{|q|^2}
\frac{|\hat\phi(q)|^2\dd q}{\lambda+|q|^2\UB_{k-1}(\lambda+|q|^2,z_{2k})}.\label{e:pranzo}
\end{equs}
We restrict the integral to the cone where $|q_2|^2\ge (1/2)|q|^2$ and we get that \eqref{e:pranzo} is lower bounded by
\begin{equ}
\frac C{f_{2k}}\int_{\mathbb R^2}
\widehat{V}(q)
\frac{\dd q}{\lambda+|q|^2\UB_{k-1}(\lambda+|q|^2,z_{2k})}
\end{equ}
where now the integral is unrestricted because the integrand depends only on $|q|$.
We can now apply again Eq. \eqref{e:ApproxDobbiamoFinire} in Lemma \ref{lemma:replacement},
so that overall~\eqref{e:pranzo} is lower bounded,
up to a multiplicative absolute constant, by
\begin{equs}[e:LBcomplete]
 \frac{1}{f_{2k+1}}&\left[\int_\lambda^1\frac{\dd \rho}{\rho\UB_{k-1}(\rho,z_{2k+1})}-\frac{\LB_{k-1}(\lambda,z_{2k+1})}{\sqrt{z_{2k+1}}}\right]
   \\
   & \geq \frac{1}{f_{2k+1}}\left[\int_\lambda^1\frac{\dd \rho}{(\rho+\rho^2)\UB_{k-1}(\rho,z_{2k+1})}-\frac{\LB_{k-1}(\lambda,z_{2k+1})}{\sqrt{z_{2k+1}}}\right]\\
&\gtrsim \frac{1}{f_{2k+1}}\left[\LB_k(\lambda,z_{2k+1})-\LB_k(1,z_{2k+1})-\frac{\LB_{k}(\lambda,z_{2k+1})}{\sqrt{z_{2k+1}}}\right]\\
&\gtrsim \frac{1}{f_{2k+1}}\left[\LB_k(\lambda,z_{2k+1})-f_{2k+1}\right]
\end{equs}
where in the second inequality we used Lemma~\ref{lemma:LBUBproperties}, and $\LB_{k-1}\leq \LB_k$,
while in the last~\eqref{e:CondK1K21} and that, for $k$ large enough, $1-1/\sqrt{z_{2k+1}}$ is bounded below by a strictly
positive constant.
Now, the $-f_{2k+1}$ just gives a constant contribution, which can be absorbed by decreasing the value of $C$
if $\lambda$ is small enough.
Using the inequality in~\eqref{e:UBFinal} for $k$ as in~\eqref{e:k}, we see that
\begin{equ}
\LB_k(\lambda,0)\gtrsim \sqrt{\Ll(\lambda,0)}\,,
\end{equ}
which, together with the definition of $f_{2k+1}$ in~\eqref{eq:K1K2},  gives
\begin{equ}
\lambda^2\tilde D(\lambda)\gtrsim (\log\Ll(\lambda,0))^{-1-\eps}\sqrt{\Ll(\lambda,0)}\,.
\end{equ}
Hence,~\eqref{eq:main2} follows at once and, by~\eqref{eq:soluzione} and the discussion thereafter,
so does Theorem~\ref{th:main}.
\end{proof}

\begin{appendix}

\section{}\label{app}

Here we collect some the technical estimates about the integrals involved in the proofs.
We also include some of the properties of the functions $\LB_k$ and $\UB_k$ from \cite[Lemma C.3]{Cannizzaro2020TheSuperdiffusivity}.
\begin{lem}\label{lemma:LBUBproperties}
For $k\in\mathbb{N}$ let $\Ll,\LB_k$ and $\UB_k$ be the functions defined in (\ref{eq:defL}) and (\ref{eq:defLBUB}).
Then, $\Ll$, $\LB_k$ and $\UB_k$ are decreasing in the first variable and increasing in the second.
For any $x>0$ and $z\geq1$, the following inequalities hold
\begin{equs}[e:bdonL]
1&\leq\LB_k(x,z)\leq\sqrt{\Ll(x,z)},\\
1&\leq\sqrt{z}\leq\sqrt{\Ll(x,z)}\leq\UB_k(x,z)\leq\Ll(x,z).
\end{equs}
Moreover for any $0<a<b$, we have
\begin{equs}
\int_a^b\frac{\dd x}{(x^2+x)\UB_k(x,z)}&=2(\LB_{k+1}(a,z)-\LB_{k+1}(b,z)),\label{e:IntUBtoLB}\\
\int_a^b\frac{\dd x}{(x^2+x)\LB_k(x,z)}&\leq2(\UB_{k}(a,z)-\UB_{k}(b,z)).\label{e:IntLBtoUB}
\end{equs}
At last, we also have
\begin{equs}[eq:LBUBderivatives]
\partial_x\Ll(x,z)&=-\frac{1}{x^2+x}\,,\quad
\partial_x\LB_k(x,z)=-\frac{1}{2(x^2+x)\UB_{k-1}(x,z)}\,,\\
\partial_x\UB_k(x,z)&=-\frac{1}{2(x^2+x)\LB_{k}(x,z)}\left(1+\frac{(\frac12\log\Ll(x,z))^k}{k!\LB_k(x,z)}\right).\,
\end{equs}
\end{lem}
\begin{proof}
All of these properties were shown in~\cite[Lemma C.3]{Cannizzaro2020TheSuperdiffusivity}. For completeness,
we add here the proof.

The two chains of inequalities in~\eqref{e:bdonL}
	are a direct consequence of the respective definitions.
	A computation of the partial derivative with respect to the second variable yields the desired monotonicity.
	Furthermore we have that
	\begin{equ}[e:Der1]
		\partial_x\Ll(x,z)=-\frac{1}{x^2+x},\,\qquad\partial_x\LB_{k}(x,z)=-\frac{1}2\frac{\LB_{k-1}(x,z)}{(x^2+x)\Ll(x,z)}
	\end{equ}
	and
	\begin{equs}[e:Der2]
		\partial_x\UB_{k}(x,z)&=-\frac{\LB_{k}(x,z)-\frac12\LB_{k-1}(x,z)}{(x^2+x)(\LB_{k}(x,z))^2}\\
		&=-\frac{1}{2(x^2+x)\LB_{k}(x,z)}\Big[1+\frac{\frac{(\frac12\log \Ll(x,z))^k}{k!}}{\LB_{k}(x,z)}\Big]\,,
	\end{equs}
	which are all strictly negative for any $x>0$ and $z\geq 1$.
	The above computation of the partial derivatives moreover reveals that
	\begin{equs}
		\int_a^b \frac{\dd x}{(x^2+x)\UB_{k}(x,z)} = 2\int_b^a \partial_x\LB_{k+1}(x,z)\dd x
		&=2\left[\LB_{k+1}(a,z)-\LB_{k+1}(b,z)\right]\,,
	\end{equs}
	which is~\eqref{e:IntUBtoLB}.
	For~\eqref{e:IntLBtoUB}, notice that
	\begin{equs}
		\int_{a}^{b}\frac{\dd x}{(x^2+x)\LB_{k}(x,z)}
		&=\int_{b}^{a}\partial_x\UB_{k}(x,z)\dd x +\frac{1}2\int_{a}^{b}\frac{ \LB_{k-1}(x,z)}{(x^2+x)\LB_{k}(x,z)^2}\dd x \\
		&\leq \int_{b}^{a}\partial_x\UB_{k}(x,z)\dd x  +\frac{1}2\int_{a}^{b} \frac{ 1}{(x^2+x)\LB_{k}(x,z)}\dd x,
	\end{equs}
	where the last inequality follows from the fact that all the terms are positive and
	for all $x$ we have $\LB_{k-1}(x,z)\leq \LB_{k}(x,z)$.
	Bringing the last term to the left hand side gives the required estimate.


\end{proof}

\begin{lem}
  \label{lemma:replacement}
  Let $V$ be a bump function satisfying Assumption~\ref{a:V}.  Let $z>1$, $f(\cdot,z)\,:\,[0,\infty)\mapsto [1,\infty)$
  be a strictly decreasing, differentiable function such that
  \begin{equ}
    \label{eq:ff'}
  -\frac{f(x)}x\le f'(x)<0 \text{ for all } x\in\R
  \end{equ}
  and $g(\cdot,z)\,:\, [0,\infty)\mapsto [1,\infty)$ a strictly decreasing function such that $g(x,z)f(x,z)\geq z$.
  Then, there exists a
  constant $C_{\Di}>0$ such that for all $z>1$, the following bound holds
  \begin{equ}\label{e:Approx}
  \left|\int_{\mathbb{R}^2}\frac{\widehat{V}(q)(\sin\theta)^2\dd q}{\lambda+|p+q|^2f(\lambda+|p+q|^2,z)}-\frac{\pi}{2}\int_{\lambda+|p|^2}^1\frac{\dd\rho}{\rho f(\rho,z)}\right|\leq
  C_{\Di}\frac{g(\lambda+|p|^2,z)}{\sqrt{z}}
\end{equ}
where $0\ne p\in\mathbb R^2$, $\theta$ is the angle between $p$ and $q$ and it is understood that the second integral is zero if $\lambda+|p|^2\ge 1$.

Moreover, for $\lambda\le 1$,
  \begin{equ}\label{e:ApproxDobbiamoFinire}
  \left|\frac12\int_{\mathbb{R}^2}\frac{\widehat{V}(q)\dd q}{\lambda+|q|^2f(\lambda+|q|^2,z)}-\frac{\pi}2\int_{\lambda}^1\frac{\dd\rho}{\rho f(\rho,z)}\right|\leq
  C_{\Di}\frac{g(\lambda,z)}{\sqrt{z}}\,.
\end{equ}

\end{lem}
\begin{proof}
  As $z$ is fixed throughout, we suppress the dependence of $f$ and $g$ on it.
  At first, we use the triangle inequality to split the left hand side of~\eqref{e:Approx} into
  \begin{equs}
  &\left|\int_{\mathbb{R}^2}\frac{\widehat{V}(q)(\sin\theta)^2\dd q}{\lambda+|p+q|^2f(\lambda+|p+q|^2)}-
  \int_{\mathbb{R}^2}\frac{\widehat{V}(q)(\sin\theta)^2\dd q}{(\lambda+|p+q|^2)f(\lambda+|p+q|^2)}\right|\label{eq:lambdabracket}\\
  &+\left|\int_{\mathbb{R}^2}\frac{\widehat{V}(q)(\sin\theta)^2\dd q}{(\lambda+|p+q|^2)f(\lambda+|p+q|^2)}-
  \int_{\mathbb{R}^2}\frac{\widehat{V}(q)(\sin\theta)^2\dd q}{(\lambda+|p|^2+|q|^2)f(\lambda+|p|^2+|q|^2)}\right|\label{eq:absolutevalues}\\
  &+\left|\int_{\mathbb{R}^2}\frac{\widehat{V}(q)(\sin\theta)^2\dd q}{(\lambda+|p|^2+|q|^2)f(\lambda+|p|^2+|q|^2)}-
  \frac{\pi}{2}\int_{\lambda+|p|^2}^1\frac{\dd\rho}{\rho f(\rho)}\right|.\label{eq:bumpfunction}
  \end{equs}
  We will bound these three terms separately.
  For the first, we re-write it as
  \begin{align*}
  \Big|\lambda\int_{\mathbb{R}^2}&\frac{\widehat{V}(q)(\sin\theta)^2(f(\lambda+|p+q|^2)-1)\dd q}{(\lambda+|p+q|^2f(\lambda+|p+q|^2))(\lambda+|p+q|^2)f(\lambda+|p+q|^2)}\Big|\\
  &\qquad\qquad\qquad\leq\lambda\int_{\mathbb{R}^2}\frac{\dd q}{(\lambda+|p+q|^2f(\lambda+|p+q|^2))(\lambda+|p+q|^2)}.
  \end{align*}
  The latter can be further split into two parts, corresponding to $|p+q|\leq|p|$ and $|p+q|>|p|$.
  In the first case $f(\lambda+|p+q|^2)\geq f(\lambda+|p|^2)$ and thus we obtain the upper bound
  \begin{align*}
  &\frac{\lambda}{f(\lambda+|p|^2)}\int_{\mathbb{R}^2}\frac{\dd q}{(\frac{\lambda}{f(\lambda+|p|^2)}+|p+q|^2)(\lambda+|p+q|^2)}\\
  &\leq\frac{\lambda}{f(\lambda+|p|^2)}
  \left(\int_{\mathbb{R}^2}\frac{\dd q}{(\frac{\lambda}{f(\lambda+|p|^2)}+|p+q|^2)^2}\right)^\frac12
  \left(\int_{\mathbb{R}^2}\frac{\dd q}{(\lambda+|p+q|^2)^2}\right)^\frac12\\
  &=\frac{C\lambda}{f(\lambda+|p|^2)}\frac{f(\lambda+|p|^2)^{\tfrac12}}{\lambda}
  =\frac{C}{\sqrt{f(\lambda+|p|^2)}}
  \leq C\frac{g(\lambda+|p|^2)}{\sqrt{z}},
  \end{align*}
  for some positive constant $C$.
  For the other case we use that $f(\lambda+|p+q|^2)
  \geq\frac{z}{g(\lambda+|p|^2)}$.
  Applying the same steps as above we get an upper bound of the form
  \begin{equ}
  \frac{C\sqrt{g(\lambda+|p|^2)}}{\sqrt{z}}
  \leq\frac{Cg(\lambda+|p|^2)}{\sqrt{z}},
  \end{equ}
 which holds as $g\geq 1$.
%

  Now we look at \eqref{eq:absolutevalues}.  First note that the
  restriction of each integral to the region $|q+p|<|p|$ has
  an upper bound of the desired form.  Indeed for the first integral
  we can use $(\sin\theta)^2\leq\frac{|p+q|^2}{|p|^2}$   (which holds for any $q_1$ and $q_3$
by elementary Euclidean geometry) to obtain
  \begin{align*}
  \int_{|p+q|<|p|}&\frac{\widehat{V}(q)(\sin\theta)^2\dd q}{(\lambda+|p+q|^2)f(\lambda+|p+q|^2)}\\
  &\leq
  |p|^{-2}\int_{|p+q|<|p|}\frac{\dd q}{f(\lambda+|p+q|^2)}
\leq  \frac{C}{f(\lambda+|p|^2)}\leq
  \frac{Cg(\lambda+|p|^2)}{\sqrt{z}}.
  \end{align*}
  For the second integral in~\eqref{eq:absolutevalues}, we can bound from above $|\sin\theta|\leq 1$,
  the denominator from below by $|p|^2f(\lambda+5|p|^2)$ and notice that the area of integration is of order $|p|^2$.

  As for the region $|q+p|\ge |p|$, define $h(x)=xf(x)$.  By \eqref{eq:ff'}, $|h'(x)|\leq 2|f(x)|$, therefore
  \begin{equs}
  |h(\lambda+|p+q|^2)-&h(\lambda+|p|^2+|q|^2)|\\
  &\leq
  2||p+q|^2-|p|^2-|q|^2|f(\min(\lambda+|p+q|^2,\lambda+|p|^2+|q|^2))\,,
  \end{equs}
  since $f$ is positive and decreasing. Therefore, we get
  \begin{align*}
  &\int_{|p+q|\geq|p|}\frac{\widehat{V}(q)(\sin\theta)^2|h(\lambda+|p+q|^2)-h(\lambda+|p|^2+|q|^2)|\dd q}{(\lambda+|p+q|^2)(\lambda+|p|^2+|q|^2)f(\lambda+|p+q|^2)f(\lambda+|p|^2+|q|^2)}\\
&\lesssim  \int_{|p+q|\geq|p|}\frac{\widehat{V}(q)(\sin\theta)^2|p||q||\cos\theta|\dd q}{(\lambda+|p+q|^2)(\lambda+|p|^2+|q|^2)f(\max(\lambda+|p+q|^2,\lambda+|p|^2+|q|^2))}\\
&\lesssim  |p|\int_{|p+q|\geq|p|}\frac{|q|\dd q}{(\lambda+|p+q|^2)(\lambda+|p|^2+|q|^2)f(\lambda+2|p|^2+2|q|^2)}\\
 &\lesssim \frac{g(\lambda+|p|^2)}{z}|p|\int_{|p+q|\geq|p|}\frac{|q|\dd q}{(\lambda+|p|^2+|q|^2)(\lambda+|p+q|^2)}\leq
    C \frac{g(\lambda+|p|^2)}{\sqrt z}
  \end{align*}
as can be seen by further splitting the last integral into the region where $|q|\ge 2|p|$ and
the complementary one, and using $z>1$.
  This concludes the estimate of the second term.

  For~\eqref{eq:bumpfunction}, we split the first integral into two regions,
  one such that $|q|^2\geq1-(\lambda+|p|^2)$ and the other given by its
  (possibly empty) complement.
  Note that on the first $\lambda+|p|^2+|q|^2\geq 1$.
  Therefore, the integral can be bounded above by
  \begin{align*}
 \frac{1}{z} \int_{|q|^2\geq1-(\lambda+|p|^2)}\widehat{V}(q)g(\lambda+|p|^2+|q|^2)\dd q\leq\frac{g(\lambda+|p|^2)}{z}\int_{\mathbb{R}^2}\widehat{V}(q)\dd q=C\frac{g(\lambda+|p|^2)}{z}\,.
  \end{align*}
  To treat the second, since $\widehat{V}(\cdot)$ is smooth and rotationally invariant,
  there is a constant $C$ such that $|\widehat{V}(q)-\widehat{V}(0)|<C|q|^2$ for $|q|\le 1$.
  We can now write the remaining integral as
  \begin{equs}
  &\int_{|q|^2<1-(\lambda+|p|^2)}\frac{(\sin\theta)^2\dd q}{(\lambda+|p|^2+|q|^2)f(\lambda+|p|^2+|q|^2)}\notag\\
  &\quad+\int_{|q|^2<1-(\lambda+|p|^2)}\frac{(\widehat{V}(0)-\widehat{V}(q))(\sin\theta)^2\dd q}{(\lambda+|p|^2+|q|^2)f(\lambda+|p|^2+|q|^2)}\,.
  \end{equs}
  By passing to polar coordinates and setting $\rho=\lambda+|p|^2+|q|^2$, the first summand can
  be immediately seen to equal the second integral in \eqref{eq:bumpfunction}.
  The second summand instead can be controlled via
  \begin{align*}
  \int_{|q|^2<1-(\lambda+|p|^2)}&\frac{|\widehat{V}(0)-\widehat{V}(q)|(\sin\theta)^2\dd q}{(\lambda+|p|^2+|q|^2)f(\lambda+|p|^2+|q|^2)}\\
  &\leq\frac{Cg(\lambda+|p|^2)}{z}\int_{|q|^2<1-(\lambda+|p|^2)}\frac{|q|^2\dd q}{\lambda+|p|^2+|q|^2}\leq C\frac{g(\lambda+|p|^2)}{z}.
  \end{align*}
  Thus, collecting all the estimates obtained so far,~\eqref{e:Approx} follows at once.

  Finally, to see \eqref{e:Approx2}, we recall that \eqref{e:Approx} holds uniformly for $p\ne0$.  Letting $p\to0$, the second integral and the r.h.s. of \eqref{e:Approx} tend to the analogous quantities in \eqref{e:Approx2}. As for the first integral in \eqref{e:Approx}, for $p\to0$ the integral over $|q|$ and $\theta$ factorizes, and we get the first integral in \eqref{e:Approx2} times $1/2$ (coming from the average of $(\sin\theta)^2$).
\end{proof}

\begin{lem}\label{lemma:off-diagonals}
Let the assumptions of Lemma \ref{lemma:replacement} be in place. 
Then, there exists a constant $C_{\oDi}>0$ such that, for every $q_1,q_2$,
\begin{equ}
|q_1|\int_{\mathbb{R}^{2}}\frac{\widehat{V}(q_3)(\sin\theta)^2\dd q_3}{(\lambda+|q_1+q_3|^2f(\lambda+|q_1+q_3|^2))|q_2+q_3|}\leq C_{\oDi}
\frac{g(\lambda+|q_1|^2)}{z}
\end{equ}
with $\theta$ the angle between $q_1$ and $q_3$.
\end{lem}
\begin{proof}
Throughout the proof the constant $C$ appearing in the bounds is independent of $q_1,\,q_2$ and $q_3$ and
might change from line to line.

We split $\mathbb{R}^2$ into three regions, $\Omega_1=\{q_3\colon|q_1+q_3|<\frac{|q_1|}{2}\}$,
$\Omega_2=\{q_3\colon|q_2+q_3|<\frac{|q_1|}{2}\}$ and
$\Omega_3=\mathbb{R}^2\setminus(\Omega_1\cup\Omega_2)$.
Note that $\Omega_1$ and $\Omega_2$ might not be disjoint,
but this is no issue as we are proving an upper bound.

In $\Omega_1$, we exploit the monotonicity of $f$ to bound
$f(\lambda+|q_1+q_3|^2)\geq f(\lambda+\frac14|q_1|^2)$.
Moreover, we estimate $(\sin \theta)^2\leq \frac{|q_3+q_1|^2}{|q_1|^2}$ and $\widehat{V}$ by a constant to get
\begin{align*}
  |q_1|\int_{\Omega_1}&\frac{\widehat{V}(q_3)(\sin\theta)^2\dd q_3}{(\lambda+|q_1+q_3|^2f(\lambda+|q_1+q_3|^2))|q_2+q_3|}\\
&\leq  C|q_1|^{-1}\int_{\Omega_1}\frac{|q_1+q_3|^2\dd q_3}{(\lambda+|q_1+q_3|^2f(\lambda+\frac14|q_1|^2))|q_2+q_3|}\\
&\leq  C\frac{|q_1|^{-1}}{f(\lambda+\frac14|q_1|^2)}\int_{\Omega_1}\frac{\dd q_3}{|q_2+q_3|}\\
&\leq\frac{C}{f(\lambda+\frac14|q_1|^2)}\leq\frac{C}{f(\lambda+|q_1|^2)}\leq C\frac{g(\lambda+|q_1|^2)}{z}\,,
\end{align*}
the last step from the third to the fourth line being a consequence of the fact that, on $\Omega_1$, $\frac{|q_1|}{2}<|q_3|<\tfrac32|q_1|$.

For $\Omega_2$ we estimate the sine differently, i.e.
\begin{equ}[e:Annoying]
    (\sin\theta)^2\leq \frac{4|q_3+q_1|^2}{|q_1|^2\vee(\frac14|q_2|^2)},
\end{equ}
which holds as, for $|q_2|\leq2|q_1|$ this is just a weaker estimate than the previous
one, while for $|q_2|\geq2|q_1|$ we claim that, in the
region $\Omega_2$ the right hand side is always greater or equal to
$1$ (and thus the inequality holds as well). Indeed, notice that since $|q_2|\geq 2|q_1|$, we have
\begin{equ}
\frac{4|q_3+q_1|^2}{|q_1|^2\vee(\frac14|q_2|^2)}=16\frac{|q_3+q_1|^2}{|q_2|^2}\,.
\end{equ}
Assume by contradiction that $|q_3+q_1|<\frac14 |q_2|$. Then
\begin{equs}
|q_3+q_2|\geq |q_2-q_1|-|q_1+q_3|>|q_2|-|q_1|-\frac14|q_2|\geq \frac{1}{4}|q_2|
\end{equs}
where in the last step we used once again that $|q_2|\geq 2|q_1|$. Now, on $\Omega_2$, $|q_2+q_3|<\frac12|q_1|$,
so that, in conclusion
\begin{equ}
 \frac{1}{4}|q_2|<|q_3+q_2|<\frac12|q_1|\leq \frac14|q_2|
\end{equ}
which is a contradiction. Hence, $|q_3+q_1|\geq\frac14 |q_2|$, from which~\eqref{e:Annoying} follows.

Plugging~\eqref{e:Annoying} into the
integral we get
\begin{align*}
    |q_1|\int_{\Omega_2}&\frac{\widehat{V}(q_3)(\sin\theta)^2\dd q_3}{(\lambda+|q_1+q_3|^2f(\lambda+|q_1+q_3|^2))|q_2+q_3|}\\
    &\leq
    C\frac{|q_1|}{|q_1|^2\vee(\frac14|q_2|^2)}\int_{\Omega_2}\frac{\dd q_3}{|q_2+q_3|f(\lambda+|q_1+q_3|^2)}.
\end{align*}
Now we use the monotonicity of $f$ to bound the previous integral from above by
\begin{equs}[e:Approx2]
    C&\frac{|q_1|}{\Big(|q_1|^2\vee(\frac14|q_2|^2)\Big)f(\lambda+(\frac32|q_1|+|q_2|)^2)}\int_{\Omega_2}\frac{dq_3}{|q_2+q_3|}\\
    &\qquad\qquad\qquad=C\frac{|q_1|^2}{\Big(|q_1|^2\vee(\frac14|q_2|^2)\Big)f(\lambda+(\frac32|q_1|+|q_2|)^2)}.
\end{equs}
We now bound this term by maximizing over $|q_2|$.
It is easy to see that it is monotonically increasing for $|q_2|<2|q_1|$.
For $|q_2|\geq2|q_1|$ we will prove that it is monotonically decreasing.
 Since $f$ satisfies assumption \eqref{eq:ff'},
 for any $a,\,b\geq 0$ we have
\begin{align*}
   & \frac{d}{dr}\left(\frac{1}{r^2f(a+(b+r)^2)}\right)=
    -\frac{2rf+2r^2(b+r)f'}{r^4f^2}\\
    &=-\frac{2}{r^3f^2}(f+r(b+r)f')\le
    -\frac{2}{r^3f}\left(1-\frac{r(b+r)}{a+(b+r)^2}\right)<0,
\end{align*}
where we suppressed the argument of $f$ and $f'$ because it does not change.
Thus, the maximum over $q_2$ of the right hand side of~\eqref{e:Approx2} is achieved at $|q_2|=2|q_1|$
and reads
\begin{equ}
    \frac{C}{f(\lambda+(\frac72|q_1|)^2)}\leq\frac{Cg(\lambda+(\frac72|q_1|)^2)}{z}\leq\frac{Cg(\lambda+|q_1|^2)}{z}\,.
\end{equ}

We are left to consider the integral over $\Omega_3$.
In this case, we first bound the $(\sin\theta)^2\le 1$ and then apply the H\"older inequality with exponents
$\frac{3}{2}$ and $3$, to the two functions $((\lambda+|q_1+q_3|^2)f(\lambda+|q_1+q_3|^2))^{-1}$ and
$|q_2+q_3|^{-1}$ with respect to the measure $\widehat{V}(q_3)\dd q_3$, so that we obtain
\begin{equs}
  &|q_1|\int_{\Omega_3}\frac{\widehat{V}(q_3)(\sin\theta
    )^2\dd q_3}{(\lambda+|q_1+q_3|^2f(\lambda+|q_1+q_3|^2))|q_2+q_3|}\notag\\
    &\leq
    |q_1|\left(\int_{\Omega_3}\frac{\widehat{V}(q_3)\dd q_3}{(\lambda+|q_1+q_3|^2f(\lambda+|q_1+q_3|^2))^\frac{3}{2}}\right)^\frac{2}{3}
    \left(\int_{\Omega_3}\frac{\widehat{V}(q_3)\dd q_3}{|q_2+q_3|^3}\right)^\frac{1}{3}\,.\label{eq:hoelder}
\end{equs}
The second integral is upper bounded by a constant factor times $|q_1|^{-1}$.

In the first integral of~\eqref{eq:hoelder}, we make the change of variables $q=q_1+q_3$,
bound the bump function $\widehat V$ by a constant and then pass to polar coordinates, hence we get
\begin{equ}\label{eq:om3OffInt}
    C\int_{\frac{|q_1|}{2}}^\infty\frac{\varrho \dd\varrho}{(\lambda+\varrho^2f(\lambda+\varrho^2))^\frac{3}{2}}\,.
\end{equ}
We split the domain of integration into two parts, $\rho^2>\lambda$ and $\rho^2\leq\lambda$ (the second one might be empty).
In the first, we note that
\begin{equ}
    \lambda+\varrho^2f(\lambda+\varrho^2)\geq \frac12(\lambda+\varrho^2)f(\lambda+\varrho^2).
\end{equ}
Using $f(x)\geq\frac{z}{g(x)}$ we obtain that this part of \eqref{eq:om3OffInt} is upper bounded by
\begin{equs}
    C\int_{\frac{\sqrt{\lambda}}{2}}^\infty\frac{\varrho \dd\varrho}{((\lambda+\varrho^2)f(\lambda+\varrho^2))^\frac{3}{2}}
    &\leq
    C\int_{\frac{|q_1|}{2}}^\infty\frac{\varrho \dd\varrho}{((\lambda+\varrho^2)f(\lambda+4\varrho^2))^\frac{3}{2}}\\&
    \leq C\left(\frac{g(\lambda+|q_1|^2)}{z}\right)^\frac{3}{2}
    \int_{\frac{|q_1|}{2}}^\infty
    \frac{\varrho \dd\varrho}
    {(\lambda+\varrho^2)^\frac{3}{2}}\\
    &\leq C|q_1|^{-1}\left(\frac{g(\lambda+|q_1|^2)}{z}\right)^{\frac32}\,,\label{eq:om3h1}
\end{equs}
where in the last step we estimated the integral by dropping $\lambda$ from the denominator.

We now turn to the second part of the integral, where $\rho^2\leq\lambda$.
We use the following
\begin{equ}
    \int_{\frac{|q_1|}{2}}^{\sqrt{\lambda}}\frac{\varrho \dd\varrho}{(\lambda+\varrho^2f(\lambda+\varrho^2))^\frac{3}{2}}\leq
    \frac{1}{f(2(\lambda+|q_1|^2))^\frac32}\int_{\frac{|q_1|}{2}}^{2\sqrt{\lambda}}\frac{ \dd\varrho}{\varrho^2}\leq C|q_1|^{-1}\frac{g(\lambda+|q_1|^2)^\frac32}{z^\frac32}\,.
\end{equ}
In conclusion, plugging these estimates into~\eqref{eq:hoelder}, we get that the integral over $\Omega_3$
is upper bounded by
\begin{equ}
    C\frac{g(\lambda+|q_1|^2)}{z}\,,
\end{equ}
and, collecting all the bounds derived so far, the statement follows at once.
\end{proof}

\begin{lem}
  \label{lemma:ff'}
The functions $\UB_k(\cdot,z)$ and $\LB_k(\cdot,z)$ satisfy the conditions of the previous lemmas.
\end{lem}
\begin{proof}
By definition $\UB_k(\cdot,z)\LB_k(\cdot,z)=\Ll(\cdot,z)\geq z$.
From Lemma \ref{lemma:LBUBproperties} we get that $\LB_k(x,z)>1$ and $\UB_k(x,z)>1$ for all $x$
and that their derivatives are both negative. Equation
\eqref{eq:LBUBderivatives} implies that
\begin{equ}
\partial_x\LB_k(x,z)=-\frac{\LB_{k-1}(x,z)}{2(x^2+x)\Ll(x,z)}\ge -\frac{\LB_k(x,z)}x
\end{equ}
and
\begin{equ}
  \partial_x\UB_k(x,z)\geq-\frac{\UB_{k}(x,z)}{(x^2+x)\Ll(x,z)}\ge -\frac{\UB_k(x,z)}x
\end{equ}
which gives \eqref{eq:ff'}.
\end{proof}

\begin{lem}\label{lemma:rhosquared}
For any $z\geq1$, $\lambda\in\mathbb{R}_+$ and $p\in\mathbb{R}^2$ such that $\lambda+|p|^2\leq 1$, we have
\begin{equ}
\left|\int_{\lambda+|p|^2}^1\frac{\dd \rho}{\rho \LB_k(\rho,z)}-\int_{\lambda+|p|^2}^1\frac{\dd \rho}{(\rho + \rho^2) \LB_k(\rho,z)}\right|\leq \frac{\UB_k(\lambda+|p|^2, z)}{z}.
\end{equ}
\end{lem}
\begin{proof}
  Note that the difference of integrals equals
  \begin{equs}
   0\le  \int_{\lambda+|p|^2}^1\frac{d\rho}{(1+\rho)\LB_k(\rho,z)}= \int_{\lambda+|p|^2}^1\frac{\UB_k(\rho,z)}{(1+\rho)\Ll(\rho,z)}d\rho\le  \frac{\UB_k(\lambda+|p|^2, z)}{z}
  \end{equs}
  because $\UB_k(\cdot,z)$ is decreasing and $\Ll(\cdot,z)\ge z$.
\end{proof}


\end{appendix}
\section*{Acknowledgements}
The authors would like to thank B\'alint T\'oth and Benedek Valk\'o for enlightening discussions,
and the anonymous referee for a very careful reading and whose comments allowed us to improve
the main statement and streamline some technical arguments.
G. C. gratefully acknowledges financial support via the EPSRC grant EP/S012524/1.
F. T. gratefully acknowledges financial support of Agence Nationale de la Recherche via the
ANR-15-CE40-0020-03 Grant LSD and of the  Austria Science Fund (FWF), Project Number P 35428-N.

\bibliographystyle{imsart-number} 
\bibliography{references}

\end{document}